\documentclass{tran-l}
\usepackage{mathrsfs}

\newtheorem{theorem}{Theorem}[section]
\newtheorem{lemma}[theorem]{Lemma}
\newtheorem{proposition}[theorem]{Proposition}
\newtheorem{corollary}[theorem]{Corollary}

\theoremstyle{definition}
\newtheorem{definition}[theorem]{Definition}

\theoremstyle{remark}

\newcount\ccount\ccount=0
\def\mkconst#1{\global\advance\ccount
    by1\edef#1{c_{\the\ccount}}\ignorespaces}
\newcount\Ccount\Ccount=0
\def\mkConst#1{\global\advance\Ccount
    by1\edef#1{C_{\the\Ccount}}\ignorespaces}

\numberwithin{equation}{section}

\begin{document}

\title{}

\author{Eric J. Olson}
\address{Department of Mathematics/084, University of Nevada, Reno, NV 89557. USA.}
\email{ejolson@unr.edu}

\author{James C. Robinson}
\address{Mathematical Institute, University of Warwick, Coventry, CV4 7AL. U.K.}
\email{jcr@maths.warwick.ac.uk}
\thanks{JCR is a Royal Society University Research Fellow, and would like to thank the Society for all their support.}


\date{March 2005}



\def\d{{\rm d}}
\def\e{{\rm e}}
\def\Re{{\mathbb R}}
\def\lap{\Delta}
\def\half{\frac{1}{2}}
\def\A{{\cal A}}
\def\be#1{\begin{equation}\label{#1}}
\def\ee{\end{equation}}
\def\bea{\begin{eqnarray*}}
\def\eea{\end{eqnarray*}}
\def\dist{{\rm dist}}
\def\etal{{\sl et al.}}
\def\<{{\langle}}
\def\>{{\rangle}}
\def\nse#1#2#3{\frac{\partial #1}{\partial t}+(#1\cdot\grad)#1-#2\lap#1+\grad
p=#3}
\def\nsefull#1#2#3{\nse{#1}{#2}{#3}\qquad\div#1=0}
\def\L{\mathbb L}
\def\N{\mathbb N}
\def\Z{\mathbb Z}
\def\grad{\nabla}
\def\div{\nabla\cdot}
\def\fa{{\rm for\ all}}
\def\AA{\mathbb A}
\def\qfa{\quad{\rm for\ all}\quad}
\def\qas{\quad{\rm as}\quad}
\def\qqor{\qquad{\rm or}\qquad}
\def\qqas{\qquad{\rm as}\qquad}
\def\qqfa{\qquad{\rm for\ all}\qquad}
\def\qqand{\qquad{\rm and}\qquad}
\def\qand{\quad{\rm and}\quad}
\def\qqwith{\qquad{\rm with}\qquad}
\def\etal{{\sl et al.}}
\def\dist{{\rm dist}}
\def\ddt#1{\frac{\d #1}{\d t}}
\def\({\left(}
\def\){\right)}
\def\cal{\mathcal}
\def\f{{\rm f}}
\def\qwith{\quad\mbox{with}\quad}
\def\qqfa{\qquad\mbox{for all}\qquad}
\def\dev{{\rm dev}}
\def\mdev{\mu\mbox{-}\dev}
\def\qqand{\qquad\mbox{and}\qquad}
\def\dB{d_{\rm B}}
\def\dA{d_{\rm A}}
\def\boxop{\mathbin{
    \hbox{\vbox{\hrule\hbox{\vrule height4pt\hskip3pt\vrule}\hrule}}}}
\def\dG{d_{\boxop}}
\def\dF{d_{\rm F}}
\def\RC{R}
\def\Qk{(I-P_k)}
\def\diam{{\rm diam}}
\def\cover{{\cal N}}
\def\xlnx#1{\Psi_{#1}}
\newenvironment{blank}[1]{}{}
\def\slog{{\rm slog}}
\def\{{\left\lbrace}
\def\}{\right\rbrace}

\title[Almost bi-Lipschitz embeddings]{Almost bi-Lipschitz
embeddings and almost homogeneous sets}

\begin{abstract}
This paper is concerned with embeddings of homogeneous spaces into
Euclidean spaces. We show that any homogeneous metric space can be
embedded into a Hilbert space using an almost bi-Lipschitz mapping
(bi-Lipschitz to within logarithmic corrections). The image of this
set is no longer homogeneous, but `almost homogeneous'. We therefore
study the problem of embedding an almost homogeneous subset $X$ of a
Hilbert space $H$ into a finite-dimensional Euclidean space. In fact
we show that if $X$ is a compact subset of a Banach space and $X-X$
is almost homogeneous then, for $N$ sufficiently large, a prevalent
set of linear maps from $X$ into $\Re^N$ are almost bi-Lipschitz
between $X$ and its image. We are then able to use the Kuratowski
embedding of $(X,d)$ into $L^\infty(X)$ to prove a similar result
for compact metric spaces.
\end{abstract}

\maketitle

\section{Introduction}
In this paper we investigate abstract embeddings between metric
spaces, Hilbert spaces, and finite-dimensional Euclidean spaces.
Historically (starting with Bouligand in 1928), attention has been
on bi-Lipschitz embeddings. By weakening this to almost bi-Lipschitz
embeddings we are able to obtain a number of new results.

A metric space $(X,d)$ is said to be $(M,s)$-{\it homogeneous\/}
(or simply {\it homogeneous\/})
if any ball of radius~$r$ can be covered
by at most $M(r/\rho)^s$ smaller balls of radius~$\rho$.
Since any subset of $\Re^N$ is homogeneous and
homogeneity is preserved under bi-Lipschitz mappings, it follows
that $(X,d)$ must be homogeneous if it is to admit a
bi-Lipschitz embedding into some $\Re^N$ (cf.~comments in Haj\l
asz, 2003). The {\it Assouad dimension} of $X$, $\dA(X)$, is
the infimum of all $s$ such that $(X,d)$ is $(M,s)$-homogeneous
for some $M\ge 1$.

\begin{blank}{
The space $(X,d)$ is homogeneous if and only if it is {\it
doubling}, i.e.~any ball of radius $r$ can be covered by at most
$K$ balls of radius $r/2$, where $K$ is independent of $r$ (see
Luukkainen, 1998).
Although, this characterisation yields
an appealingly simple definition,
we will need to base our generalised notion of almost
homogeneous spaces on the original definition of homogeneity.
}\end{blank}

Assouad (1983) showed that $(X,d)$ is homogeneous if and only if the
snowflake spaces $(X,d^\alpha)$ with $0<\alpha<1$ admit bi-Lipschitz
embeddings into some $\Re^N$ (where $N$ depends on $\alpha$).
However, there are examples due to Laakso (2002; see also Lang \&
Plaut, 2001) of homogeneous spaces that do not admit a bi-Lipschitz
embedding into any $\Re^N,$ nor even into an infinite-dimensional
Hilbert space. This paper starts with a simple result, based on
Assouad's argument, that any homogeneous metric space admits an {\it
almost bi-Lipschitz} embedding into an infinite-dimensional Hilbert
space.

The class of $\gamma$-almost $L$-bi-Lipschitz mappings
$f:(X,d)\rightarrow(\tilde X,\tilde d)$ (or {\it almost
bi-Lipschitz mappings\/} for short) consists of all those maps for
which there exists a $\gamma\ge0$ and an $L>0$ such that
\be{abbi}
\frac{1}{L}\,\frac{d(x,y)}{\slog(d(x,y))^\gamma}
\le \tilde d(f(x),f(y))\le L\, d(x,y)
\ee
for all $x,y\in X$ such that $x\ne y$. Here $\slog(x)$ is the
`symmetric logarithm' of $x$, defined as
$$
\slog(x):=\log(x+x^{-1}),
$$
and so an almost bi-Lipschitz map is bi-Lipschitz to within
logarithmic corrections.

\blank{Here $\xlnx{\gamma}:[0,\infty)\to[0,\infty)$ is a
continuous monotone increasing function such that $x\le
\xlnx{\gamma}(x)$, $\xlnx{\gamma}(x)\sim x(-\log x)^{\gamma}$ as
$x\to 0$ and $\xlnx{\gamma}(x)\sim x(\log x)^{\gamma}$ as $x\to
\infty.$ Note that the function
$\Psi_\gamma(x)=x\big[\log(e+x+1/x)\big]^\gamma$ is a example of
such a function.}

Although the bi-Lipschitz image of a homogeneous set is
homogeneous, this is not true for almost bi-Lipschitz images; they
are, however, {\it almost homogeneous}: we say that $(X,d)$ is
$(\alpha,\beta)$-almost $(M,s)$-homogeneous if
\be{aHom}
\cover_X(r,\rho)\le M
    \left(\frac{r}{\rho}\right)^s\,\slog(r)^\beta\slog(\rho)^\alpha
\ee
for all $0<\rho<r<\infty.$
The {\it Assouad $(\alpha,\beta)$-dimension\/} of $X$,
$\dA^{\alpha,\beta}(X)$, is the infimum of all $s$ such that $X$ is
$(\alpha,\beta)$-almost $(M,s)$-homogeneous for some $M\ge 1$.

Olson (2002) showed that given a compact $X\subset\Re^N$ with
$d_{\rm A}(X-X)=d$ then almost every projection of rank~$k>d$
provides an almost bi-Lipschitz embedding of $X$ into $\Re^k$. In
this paper we show a similar result for compact subsets $X$ of a
Hilbert space: If the set of differences\footnote{The introduction
of a condition on the dimension of the set $X-X$ of differences,
rather than on $X$ itself, is common in the literature on abstract
embeddings. The proof of Ma\~n\'e's 1981 embedding theorem requires
the Hausdorff dimension of $X-X$ to be finite, a condition not
ensured by the finiteness of $d_{\rm H}(X)$. Foias \& Olson (1996)
and Hunt \& Kaloshin (1999) treat the upper box-counting dimension
which is unusual in having the property that $\dF(X)<\infty$ implies
that $\dF(X-X)<\infty$. [Recall that
$\dF(X)=\limsup_{\epsilon\rightarrow0}\log
N(X,\epsilon)/(-\log\epsilon)$, where $N(X,\epsilon)$ is the minimum
number of balls of radius $\epsilon$ needed to cover $X$.]} $X-X$ is
almost homogeneous with $\dA^{\alpha,\beta}(X-X)=d$ then `most'
linear maps into Euclidean spaces $\Re^k$ of with $k>d$ provide
almost bi-Lipschitz embeddings of $X$.  More explicitly, if $k>d$
then the set of almost bi-Lipschitz embeddings into $\Re^k$ is
prevalent in the space of all linear maps into $\Re^k$, in the sense
of Hunt, Sauer \& Yorke (1992). We then extend this result to
subsets of Banach spaces.

There is an unfortunate gap here. An almost homogeneous metric
space has an almost bi-Lipschitz image that is an almost
homogeneous subset of a Hilbert space. However, our embedding
theorem for a subset $X$ of a Hilbert space requires that not $X$
itself, but the set $X-X$ of differences  is almost homogeneous.

By using the Kuratowski isometric embedding of $(X,d)$ into
$L^\infty(X)$ we can assign a meaning to ``$X-X$" even when $X$ is a
metric space. With this interpretation, we can also show that if
$(X,d)$ is a compact metric space then the assumption that $X-X$ is
almost homogeneous is sufficient to ensure that $(X,d)$ can be
embedded into a Euclidean space in an almost bi-Lipschitz way.

In Section 2 we state some elementary properties of the
$(\alpha,\beta)$-Assouad dimension and show that any almost
homogeneous metric space $(X,d)$ can be embedded into a Hilbert
space in an almost bi-Lipschitz way; that such almost bi-Lipschitz
images of almost homogeneous spaces are again almost homogeneous is
shown in Section 3. Section 4 treats the local versions of
homogeneity and almost homogeneity. Section 5 contains our main
result on embedding a subset $X$ of a Hilbert space with $X-X$
almost homogeneous, while in Section 6 we consider what is possible
for such subsets knowing only properties of $X$. In Section 7 we
generalise our main theorem to treat subsets of Banach spaces, and
hence give a result for compact metric spaces. In Section 8 we
explore the relationship between $\dA^{\alpha,\beta}(X)$ and
$\dA^{\alpha,\beta}(X-X)$. After Section 9, where we give an example
of a homogeneous set that cannot be bi-Lipschitz embedded into any
$\Re^k$ using any {\it linear} map, we finish with some interesting
open problems.

\section{Almost homogeneous metric spaces}\label{sec_AH}

As discussed above, we will say that a metric space $(X,d)$ is
$(\alpha,\beta)$-almost $(M,s)$-homogeneous (or simply almost
homogeneous) if any ball of radius $r$ can be covered by at
most\footnote{For bounded metric spaces (\ref{axHom}) could be
replaced by
$$
\cover_X(r,\rho)\le
M'\left(\frac{r}{\rho}\)^s\,\log(\e+\rho^{-1})^{\gamma},
$$
(in terms of our current definition we would have $M'\ge M$ and
$\gamma=\alpha+\beta$) while for compact spaces the factor of $\e$
in the logarithm could also be dropped by considering only
$\rho\le r\le\epsilon$ for some $\epsilon>0$ (see Section
\ref{local}). However, (\ref{axHom}) allows us to treat general
metric spaces.}
\be{axHom}
 \cover_X(r,\rho)\le M\left(\frac{r}{\rho}\)^s\,\slog(r)^\beta\,\slog(\rho)^\alpha
\ee
balls of radius $\rho$ (with $\rho<r$), for some $M\ge1$ and
$s\ge0$, where $\slog(x)=\log(x+x^{-1})$.

We now give some simple properties of the function slog.

\begin{lemma}\label{SlogLem}
Given $L>0$ and $\gamma\ge 0$, there exist constants $A_L, B_L,
a_\gamma, b_\gamma, \sigma\in (0,\infty)$ independent of $x$ such
that
\begin{enumerate}
\item[(p1)]
$|\log(x)|\le\slog(x)\le\log 2+|\log(x)|$, in particular
$\slog(2^k)\le(1+|k|)\log 2$,
\item[(p2)]
$A_L\,\slog(x)\le \slog(L x)\le B_L\,\slog(x)$,
\item[(p3)]
$a_\gamma\,\slog(x)\le\slog(x\,\slog(x)^\gamma)\le
b_\gamma\,\slog(x)$,
\end{enumerate}
for all $x\ge 0$, and
\begin{enumerate}
\item[(p4)] if $2^{-(k+1)}\le x\le 2^{-k}$ then
$\slog(x)\ge\sigma\,\slog(2^{-k})$.
\end{enumerate}
\end{lemma}

\begin{proof}
(p1) is elementary. For (p2) consider the quotient function
$g:(0,\infty)\to (0,\infty)$ defined by
$$
    g(x)=\frac{\slog(Lx)}{\slog(x)}.
$$
Let $a_L=\inf\{\, g(x) : x\in (0,\infty)\,\}$ and $b_L=\sup \{\,
g(x) : x\in (0,\infty)\,\}$. Since
$$\lim_{x\to 0}g(x)=1,\quad\lim_{x\to\infty}g(x)=1,
\quad\hbox{and}\quad 0<g(x)<\infty\hbox{ for }x\in (0,\infty),$$
then both $a_L$ and $b_L$ are finite positive constants. The proof
of (p3) is similar. For (p4)
set $x=2^{-r}$ with $k\le r\le k+1$.
Since $\slog(x)=\log(x+1/x)\ge\log 2$
and $\slog(2^{-r})\ge |\log 2^{-r}|=|r|\log 2$ from (p1),
then
$\slog(x)\ge (1+|r|)/2$.
Therefore, the estimate
$$
\frac{\slog(2^{-k})}{\slog(x)}
    \le
    \frac{(1+|k|)\log 2}{(1+|r|)/2}\le 4\log 2
$$
gives (p4) with $\sigma=1/(4\log2)$.
\end{proof}

We define the Assouad $(\alpha,\beta)$-dimension of $X$, $d_{\rm
A}^{\alpha,\beta}(X)$, to be the infimum of all $s$ for which $X$ is
$(\alpha,\beta)$-almost $(M,s)$-homogeneous. When $\alpha=\beta=0$
we recover the standard definition of a homogeneous space and the
usual Assouad dimension.

We note here that it is straightforward to show that the Assouad
$(\alpha,\beta)$-dimension satisfies the minimal properties we
would ask for in a dimension, namely that
$$
X\subseteq Y\quad\Rightarrow\quad d_{\rm A}^{\alpha,\beta}(X)\le
d_{\rm A}^{\alpha,\beta}(Y),\qquad d_{\rm A}^{\alpha,\beta}(X\cup
Y)=\max(d_{\rm A}^{\alpha,\beta}(X),d_{\rm A}^{\alpha,\beta}(Y)),
$$
and $d_{\rm A}^{\alpha,\beta}({\cal O})=n$ if ${\cal O}$ is an open
subset of $\Re^n$. Furthermore,
\be{ordered}
\alpha_1\ge \alpha_2\qand \beta_1\ge\beta_2\quad\Rightarrow\quad
\dA^{\alpha_1,\beta_1}(X)\le \dA^{\alpha_2,\beta_2}(X).
\ee

\medskip

We now show that if $(X,d)$ is almost homogeneous then
it can be embedded into an infinite-dimensional Hilbert space in
an almost bi-Lipschitz way. Key to this result is the following
proposition, which although not given explicitly in this form,
essentially occurs in Assouad's paper. Indeed, it is the main
ingredient in his proof of the existence of bi-Lipschitz maps
between $(X,d^\alpha)$ and $\Re^N$.

\begin{proposition}\label{assouad}
Let $(X,d)$ be an $(\alpha,\beta)$-almost $(M,s)$-homogeneous
metric space and distinguish a point $a\in X$. Then there are
constants $A,B,C>0$ such that for every $j\in\Z$ there exists a
map $\phi_j:(X,d)\rightarrow\Re^{M_j}$, where $M_j=
C(1+|j|)^{\alpha+\beta} $, with $\phi_j(a)=0$, and for every
$x_1,x_2\in X$
\begin{enumerate}
\item[(a1)] $2^{-(j+1)}<d(x_1,x_2)\le
2^{-j}$ implies that $\|\phi_j(x_1)-\phi_j(x_2)\|\ge A$, and
\item[(a2)]  $\|\phi_j(x_1)-\phi_j(x_2)\|\le
BM_j \min[1,2^jd(x_1,x_2)]$.
\end{enumerate}
\end{proposition}

\begin{proof}
The proof follows exactly the steps in Assouad's original paper
(see also the lecture notes of Heinonen (2003) for an account that
is easier to follow) which we outline very briefly here: if $N_j$
is a maximal $2^{-j}$ net in $(X,d)$, then for every $x\in X$
\bea
{\rm card}\left(N_j\cap B(x,12\cdot 2^{-j})\right)
&\le& \cover_X(12\cdot 2^{-j},2^{-j-1}) \\
&\le& 24M\,\slog(12\cdot
2^{-j})^\alpha\,\slog(2^{-j-1})^\beta\\
&\le& C(1+|j|)^{\alpha+\beta}
\eea
where the constant $C$ is a product of $M$ and the constants
appearing in Proposition \ref{SlogLem}. Thus, there exists a
`colouring map' $\kappa_j:N_j\rightarrow\{e_1,\ldots,e_{M_j}\}$,
where $\{e_1,\ldots,e_{M_j}\}$ is the standard basis of
$\Re^{M_j}$, such that $\kappa_j(a)\neq\kappa_j(b)$ if
$d(a,b)<12\cdot 2^{-j}$. Let
$$
\tilde \phi_{j}(x)=\sum_{a_i\in
N_j}\max\{(2-2^j d(x,a_i)),0\}\kappa_j(a_i).
$$
Note that $2^{2-j}<d(x_1,x_2)\le 2^{3-j}$ implies
$\tilde \phi_j(x_1)$ is orthogonal to $\tilde \phi_j(x_2)$.
It is then straightforward to show that the map
$\phi_j(x)=\tilde\phi_{j+3}(x)-\tilde\phi_{j+3}(a)$
satisfies the properties given in the statement of the
proposition.
\end{proof}

\begin{theorem}\label{agone}
Let $(X,d)$ be an $(\alpha,\beta)$-almost $(M,s)$-homogeneous
metric space and $H$ an infinite-dimensional separable Hilbert
space. Then, for every $\gamma>\alpha+\beta+\half$, there exists a
map $f:X\to H$ and a constant $L$ such that
$$
\frac{1}{L}\,\frac{d(x,y)}{\slog(d(x,y))^\gamma}
    \le\|f(s)-f(t)\|\le L\,d(x,y),
$$
i.e.,~$f$ is $\gamma$-almost bi-Lipschitz.
\end{theorem}

\begin{proof}
Let $\{e_j\}_{j\in\Z}$ be an orthonormal set of vectors in some
Hilbert space. Let $\delta>1/2$ and define
$f:(X,d)\rightarrow\bigoplus_{j=1}^\infty \Re^{M_j}\otimes
e_{j}\simeq H$ by
\be{series}
f(x)=\sum_{j=-\infty}^\infty \frac{2^{-j}}{(1+|j|)^{\delta}M_j}\,
\phi_j(x)\otimes e_j,
\ee
where the maps $\phi_j$ are those of Proposition \ref{assouad}.
Since $f(a)=0$, then the upper bound on $\|f(s)-f(t)\|$ that
we now prove will also show convergence of the series~(\ref{series})
defining $f$.
Let $(x_1,x_2)$ be a pair of distinct points of $X$.
Thus,
there exists $l\in \Z$ such that
$2^{-(l+1)}<d(x_1,x_2)\le2^{-l}$. Note that for such a pair of points
$\|\phi_l(x_1)-\phi_l(x_2)\|\ge A$.
We have
\mkconst\Cagone
\bea
\|f(x_1)-f(x_2)\|^2&
    =&\sum_{j=-\infty}^\infty
    \frac{2^{-2j}}{(1+|j|)^{2\delta}}\,
        \frac{\|\phi_j(x_1)-\phi_j(x_2)\|^2}{M_j^2} \\
&\le&\sum_{j=-\infty}^\infty \frac{B^2}{(1+|j|)^{2\delta}}\,
    d(x_1,x_2)^2\\
&\le& \Cagone\, d(x_1,x_2)^2,
\eea
where the sum converges since $2\delta>1$.

\mkconst\Cagtwo\mkconst\Cagthree

The lower bound is straightforward, since
\bea
\|f(x_1)-f(x_2)\|
    &\ge&\frac{2^{-l}}{(1+|l|)^\delta M_l}\,
    \|\phi_l(x_1)-\phi_l(x_2)\|\ \ge\ A\, \frac{2^{-l}}{(1+|l|)^\delta M_l}  \\
    &\ge&\Cagtwo\,
        \frac{2^{-l}}{(1+|l|)^{\alpha+\beta+\delta}}
   \ \ge\
   \Cagtwo\frac{d(x,y)}{(1+|l|)^{\alpha+\beta+\delta}}\\
\eea
Since $d(x,y)=2^{-r}$ with $l\le r<l+1$ it follows using (p1) from
Lemma \ref{SlogLem} that
$$
\frac{1+|l|}{\slog(d(x,y))}=\frac{1+|l|}{\slog (2^{-r})}\ge
\frac{1+|l|}{(1+|r|)\log 2}\ge\frac{1}{2\log 2},
$$
and so
$$
\|f(x_1)-f(x_2)\|\ge
\Cagthree\frac{d(x,y)}{\slog(d(x,y))^{\alpha+\beta+\delta}}.
$$
Taking $L=\max(\Cagone,1/\Cagthree)$ finishes the proof.
\end{proof}

We note here that if $(X,d)$ is bounded then there exists a $k$ such
that $d(x_1,x_2)\le 2^{k}$ for all $x_1,x_2\in X$.
In this case the definition of $f$ in \eqref{series} can be
simplified to
\be{seriescomp}
f(x)=\sum_{j=-k}^\infty \frac{2^{-j}}{(1+|j|)^\delta M_j}\,
\phi_j(x)\otimes e_j
\ee
and will still provide a $\gamma$-almost bi-Lipschitz embedding.

%

\section{Almost bi-Lipschitz images of sets}

\def\f{{\rm f}}

Since we can embed any almost homogeneous metric space into a
Hilbert space using an almost bi-Lipschitz map, it is natural to
study the effect of such mappings on almost homogeneous spaces.
Here we show that almost bi-Lipschitz images of almost homogeneous
metric spaces are still almost homogeneous. In particular this
implies that it is necessary that $X$ be almost homogeneous if it
is to enjoy an almost bi-Lipschitz embedding into some $\Re^N$.

\begin{lemma}\label{biLip}
Let $(X,d)$ be an $(\alpha,\beta)$-almost $(M,s)$-homogeneous
metric space
and $\phi:(X,d)\rightarrow(\tilde X,\tilde d)$ a $\gamma$-almost
$L$-bi-Lipschitz map. Then $(\phi(X),\tilde d)$ is an almost
homogeneous metric space with $
\dA^{\alpha+\gamma,\beta+\gamma}(X)\le
\dA^{\alpha,\beta+\gamma}(\phi(X))\le \dA^{\alpha,\beta}(X)$.
\end{lemma}


\begin{proof}
\ccount=0 Increase $L$ if necessary so that
\be{Lup}
L^2b^\gamma(\log 2)^\gamma\ge 1,
\ee
where here and in the rest of the proof $b=b_\gamma$, where
$b_\gamma$ is the constant occurring in (p3) in Lemma \ref{SlogLem};
clearly $\phi$ remains $\gamma$-almost $L$-bi-Lipschitz
under this assumption.
Take $s>\dA^{\alpha,\beta}(X)$, $0<\rho<r<\infty$, and
consider an arbitrary ball $B_{\tilde X}(\phi(x),r)$ of radius $r$
in $\phi(X)$. Now, we have\mkconst\CbiLip
$$
B_{\tilde X}(\phi(x),r)\subseteq
\phi\{B_X(x,Lrb^\gamma\slog(Lrb^\gamma)^\gamma)\},
$$
since using (p3) in Lemma \ref{SlogLem}
$$
\frac{1}{L}\frac{Lrb^\gamma\slog(Lrb^\gamma)^\gamma}{\slog(Lrb^\gamma\slog(Lrb^\gamma)^\gamma)^\gamma}
\ge\frac{rb^\gamma\slog(Lrb^\gamma)^\gamma}{[b\,\slog(Lrb^\gamma)]^\gamma}=r.
$$
By our choice of $L$ in (\ref{Lup}) and since $\rho<r$ we have
$0<\rho/L<Lrb^\gamma\slog(Lrb^\gamma)^\gamma$ and so we can cover
$B_X(x,Lrb^\gamma\slog(Lrb^\gamma)^\gamma)$ by
\bea
&&    \cover_X(Lrb^\gamma\slog(Lrb^\gamma)^\gamma,\rho/L)\\
&&\qquad\qquad\le M
        \(\frac{Lrb^\gamma\,\slog(Lrb^\gamma)^\gamma}{\rho/L}\)^s\,\slog(Lrb^\gamma\,\slog(Lrb^\gamma)^\gamma)^\beta\,\slog(\rho/L)^\alpha\\
&&\qquad\qquad \le
\CbiLip\(\frac{r}{\rho}\)^s\,\slog(r)^\beta\slog(\rho)^\alpha
\eea
balls of radius $\rho/L$ (in $X$) where $\CbiLip$ depends on $M$,
$L$ and the constants appearing in Lemma \ref{SlogLem}. Denote
these balls by $B_X(x_i,\rho/L)$. Since
$$
    \phi\{B_X(x_i,\rho/L)\}\subseteq B_{\tilde X}(\phi(x_i),\rho)$$
and $B_{\tilde X}(\phi(x),r)$ was arbitrary, it follows that
$$
    \cover_{\phi(X)}(r,\rho)\le \CbiLip\(\frac{r}{\rho}\)^s\,\slog(r)^\beta\slog(\rho)^\alpha
$$
for any $0<\rho<r<\infty$. Thus $\phi(X)$ is
$(\alpha,\beta+\gamma)$-almost $(\CbiLip,s)$-homogeneous. Taking
the infimum over $s>\dA^{\alpha,\beta}(X)$ yields
$\dA^{\alpha,\beta+\gamma}(\phi(X))\le \dA^{\alpha,\beta}(X)$.

By considering similarly the inverse map $\phi^{-1}:\phi(X)\to X$
one obtains the lower bound
$\dA^{\alpha,\beta+\gamma}(\phi(X))\ge\dA^{\alpha+\gamma,\beta+\gamma}(X)$.
\end{proof}

Combined with Lemma \ref{biLip}, the embedding result of
Proposition \ref{agone} shows that any almost homogeneous metric
space $(X,d)$ has an almost bi-Lipschitz image $f(X)$ that is an
almost homogeneous subset of a Hilbert space.


We end
by noting since almost bi-Lipschitz maps are, in fact, Lipschitz
then for any almost bi-Lipschitz map $\phi$ the upper box-counting
(`fractal') dimension satisfies $\dF(\phi(X))\le \dF(X)$.
Moreover, it is not difficult to prove the following:

\begin{lemma}
Let $(X,d)$ be a metric space
and $\phi:(X,d)\rightarrow(\tilde X,\tilde d)$ an
almost bi-Lipschitz map. Then
$\dF(\phi(X))=\dF(X)$.
\end{lemma}

\section{Aside: Compact spaces and local versions of (almost) homogeneity}\label{local}

In this section we briefly discuss the local definitions of
homogeneity and almost homogeneity, and the dimensions associated
with them. While they agree for compact spaces, they are distinct
in general.

A metric space $(X,d)$ is said to be {\it locally
$(M,s)$-homogeneous\/} (or simply {\it locally homogeneous\/}) if
there exists an $\epsilon>0$ such that any ball of radius
$r<\epsilon$ can by covered by at most $M(r/\rho)^s$ smaller balls
of radius $\rho$. The introduction of the constant $\epsilon$ for a
locally homogeneous space may be interpreted as the small scale
beneath which the set may be viewed as homogeneous. In this case $M$
may depend on $\epsilon$ which in turn depends on the units of
measurement used in the definition of the metric.

Movahedi-Lankarani (1992) defined the metric (or `Bouligand')
dimension
\be{DEFB}
\dB(X)=\lim_{\epsilon\to 0}\lim_{t\to \infty}
    \sup\{
        \frac{\log \cover_X(r,\rho)}{\log (r/\rho)}
        :0<\rho<r<\epsilon\hbox{ and } r> t\rho\},
\ee
where $\cover_X(r,\rho)$ is the minimum number of balls of radius
$\rho$ necessary to cover any ball of radius $r$. This dimension,
$\dB(X)$, is the infimum of all $s$ such that $(X,d)$ is locally
$(M,s)$-homogeneous for some $M\ge 1$.

Here we give a simple example that shows that the concepts of
homogeneous and locally homogeneous are indeed different. Let $H$
be a Hilbert space with orthonormal basis given by
$\{e_n\}_{n\in\N}$. Define
$$
    X=\{\, \rho_n e_n : n\in\N\,\}
\quad\hbox{where}\quad \rho_n=1-\frac{1}{n}.
$$
If $(X,d)$ is $(M,s)$-homogeneous for some $M$ and $s$ then
\be{FORCONT}
    \cover_X(\rho_{2n},\rho_n)\le M(\rho_{2n}/\rho_n)^s
        =M\(\frac{2n-1}{2n-2}\)^s \le M.
\ee
However, each ball $B(0,\rho_{2n})$ contains the $n$ points
$$\{0\}\cup \{\, \rho_k e_k : n<k<2n \,\}$$
which are mutually more than a distance $\rho_n$ apart. Therefore
$\cover_X(\rho_{2n},\rho_n)\ge n$. Taking $n$ large enough shows
that \eqref{FORCONT} cannot hold, and so $(X,d)$ is not
homogeneous. On the other hand, $(X,d)$ is locally homogeneous for
any $\epsilon<1$.

Note that if $(X,d)$ is compact, then the notions of homogeneous
and locally homogeneous are equivalent (see Olson, 2002). Thus
$\dA(X)=\dB(X)$ for compact spaces $X$.

As with homogeneous spaces, there is a similarly distinct notion
of {\it locally $(\alpha,\beta)$-almost $(M,s)$-homogeneous}.
This means there is some $\epsilon>0$ such that \eqref{axHom} holds
for all $0<\rho<r<\epsilon$. Similar arguments to those given in
Olson (2002) show that the notions of almost homogeneous and
locally almost homogeneous are equivalent when $(X,d)$ is compact.
Define the {\it local Assouad $(\alpha,\beta)$-dimension\/} of
$X$, $\dB^{\alpha,\beta}(X),$ to be the infimum of all $s$ such
that $(X,d)$ is locally $(\alpha,\beta)$-almost
$(M,s)$-homogeneous for some $\epsilon>0$ and $M\ge 1$.

\blank{We now state a few elementary inequalities involving the
almost linear functions $\Psi_\gamma(x)$ described in the
introduction.

\begin{proposition}\label{PsiProp}
Given $s,L>0$ and $\gamma,\alpha,\beta\ge 0$, then there exists
constants $a,b\in(0,\infty)$ independent of $x$ such that
\begin{enumerate}
\item[(p1)]
$a\Psi_\gamma(x)\le \Psi_\gamma(L x)\le b\Psi_\gamma(x),$
\item[(p2)]
$a x^{s-1}\Psi_{\gamma s}(x) \le \big[\Psi_\gamma(x)\big]^s\le b
x^{s-1}\Psi_{\gamma s}(x),$
\item[(p3)]
$a x^2 \le \Psi_\gamma(x)\Psi_{\gamma}^{-1}(x) \le b x^2$
\item[(p4)]
$a \Psi_{\alpha+\beta}(x)\le (\Psi_{\alpha}\circ\Psi_{\beta})(x)
\le b \Psi_{\alpha+\beta}(x)$, and
\item[(p5)]
$a x \Psi_{\alpha+\beta}(x)\le \Psi_{\alpha}(x)\Psi_{\beta}(x) \le
b x \Psi_{\alpha+\beta}(x)$
\end{enumerate}
for all $x\ge 0$.
\end{proposition}

\begin{proof}
Each part is similar.  We prove (p1) as an example. Consider the
quotient function $g:(0,\infty)\to (0,\infty)$ defined by
$$
    g(x)=\frac{\Psi_\gamma(Lx)}{\Psi_\gamma(x)}.
$$
Let $a=\inf\{\, g(x) : x\in (0,\infty)\,\}$ and $b=\sup \{\, g(x)
: x\in (0,\infty)\,\}$. Since
$$\lim_{x\to 0}g(x)=L,\quad\lim_{x\to\infty}g(x)=L,
\quad\hbox{and}\quad 0<g(x)<\infty\hbox{ for }x\in (0,\infty),$$
then both $a$ and $b$ are finite positive constants. The proof for
each of the other statements is nearly identical.
\end{proof}}

Let $(X,d)$ be a metric space. In general $\dB^{\alpha,\beta}(X)\le
\dA^{\alpha,\beta}(X)$.
Both $\dA^{\alpha,\beta}$ and $\dB^{\alpha,\beta}$ are invariant
under a rescaling of the metric. Thus, the metric space $(\tilde
X,\tilde d)$ where $\tilde X=X$ and $\tilde d= \eta d$ for some
$\eta>0$ has $\dA^{\alpha,\beta}(\tilde X)=\dA^{\alpha,\beta}(X)$
and $\dB^{\alpha,\beta}(\tilde X)=\dB^{\alpha,\beta}(X)$. Note that
$$
    \dB^{\alpha+\theta\beta,(1-\theta)\beta}(X)
    \le \dB^{\alpha,\beta}(X)\\
    \le \dB^{(1-\theta)\alpha,\theta\alpha+\beta}(X)
$$
for $0\le \theta\le 1$. Moreover, if $X$ is compact, then
$$  \dF(X)\le
    \dA^{\alpha,\beta}(X)=\dB^{\alpha,\beta}(X)$$
where $\dF(X)$ denotes the fractal or upper box-counting dimension.

We note here that $\dB$ shares with $\dA$ the usual properties of
dimension discussed in Section \ref{sec_AH}, along with the
monotonicty property in (\ref{ordered}).

\section{Embedding Hilbert subsets $X$ with $X-X$ homogeneous}

In this section we prove our main result, in which we take a
subspace $X$ of a Hilbert space, assume that $X-X$ is almost
homogeneous, and obtain an almost bi-Lipschitz embedding into a
finite-dimensional space.

Our argument is essentially a combination of that of Olson (2002),
who treated a subset $X$ of a Euclidean space with $\dA(X-X)$
finite, and that of Hunt \& Kaloshin (1999), who considered a
subset of a Hilbert space with finite upper box-counting
(`fractal') dimension. The key to combining these successfully is
Lemma \ref{stlem}, below.

In line with the treatment in Sauer et al.~(1991) and in Hunt \&
Kaloshin (1999), our main theorem is expressed in terms of
prevalence. This concept, which generalises the notion of `almost
every' from finite to infinite-dimensional spaces, was introduced by
Hunt, Sauer \& Yorke (1992); see their paper for a detailed
discussion.

\begin{definition}
A Borel subset $S$ of a normed linear space $V$ is
\emph{prevalent} if there exists a compactly supported
probability measure $\mu$ such that $\mu(S+v)=1$
for all $v\in V$.
In particular, if $S$ is prevalent then $S$ is dense in $V$.
\end{definition}

Note that if we
set $Q={\rm supp}(\mu)$ then $Q$ can be thought of as a `probe
set', which consists of `allowable perturbations' with which,
given a $v\in V$, we `probe' and test whether $v+q\in S$ for
almost every $q\in Q$.

Since we will use it below, and for its historical importance, we
quote Hunt \& Kaloshin's result here, in a form suitable for what
follows. Given a set $X$, we recall here that its upper
box-counting (`fractal') dimension is defined as
$$
\dF(X)=\limsup_{\epsilon\rightarrow0}\frac{\log
N(X,\epsilon)}{-\log\epsilon},
$$
where $N(X,\epsilon)$ denotes the minimum number of balls of
radius $\epsilon$ necessary to cover $X$; and its thickness
exponent, $\tau(X)$, is
\be{thickness}
\tau(X)=\limsup_{\epsilon\rightarrow0}
    \frac{\log d(X,\epsilon)}{-\log\epsilon},
\ee
where $d(X,\epsilon)$ is the minimum dimension of all
finite-dimensional subspaces, $V$, of $B$ such that every point of
$X$ lies within $\epsilon$ of $V$. We note here for later use that
$\tau(X)\le\dF(X)$.

\begin{theorem}[Hunt \& Kaloshin]\label{HK}  Let $X$ be a compact
subset of a Hilbert space $H$, $D$ an integer with $D>\dF(X-X)$, and
$\tau(X)$ the thickness exponent of $X$. If $\theta$ is chosen with
$$\theta>\frac{D(1+\tau(X)/2)}{D-\dF(X-X)}$$
then for a prevalent set of linear maps $L:B\rightarrow\Re^D$ there
exists a $c>0$ such that
$$
c\|x-y\|^\theta\le|Lx-Ly|\le \|L\|\|x-y\|\quad\mbox{for all}\quad x,y\in
X;
$$
in particular these maps are injective on $X$.
\end{theorem}

We note here that $\dF(X-X)\le2\dF(X)$, so that for zero thickness
sets with finite box-counting dimension one can choose any
$D>2\dF(X)$ and
$\theta>D/(D-2\dF(X))$.

\subsection{Construction of the probability measure $\mu$ for a
given $X$}\label{measure}

We now apply the definition of prevalence given a particular
compact subset $X$ of our Hilbert space $H$ such that $X-X$ is
$(\alpha,\beta)$-almost $(M,s)$-homogeneous.

For some fixed $N$, let $V$ be the set of linear functions $L:H\to
\Re^N$. We now construct a compactly supported probability measure
$\mu$ on $V$ (as required by the definition of prevalence) that is
carefully tailored to the particular set $X$. Key to this is the
following result.

\mkconst{\lemC}
\begin{lemma}\label{stlem}
Suppose that $X$ is a compact $(\alpha,\beta)$-almost
$(M,s)$-homogeneous subset of $H$. Then there exists a sequence of
nested linear subspaces $U_n$ with $U_{n}\subseteq U_{n+1}$,
$$\dim U_n\le C(1+n)^{\alpha+\beta+1},$$
and
$$
\|P_nx\|\ge\frac{1}{8}\|x\|\mbox{\ \ for all } x\in X\
\mbox{with}\ \|x\|\ge 2^{-n},
$$
where $P_n$ is the orthogonal projection onto $U_n$.
\begin{blank}{
 Furthermore
if $U=\cup_{k=1}^\infty U_k$ and $P$ is the orthogonal projection
onto $U$ then
$$
|Px|\ge\frac{1}{8}|x|\mbox{\ \ for all}\ x\in X.
$$
}\end{blank}
\end{lemma}

\begin{proof}
Consider the collection of shells
$$
\Delta_j=\{x\in X:2^{-(j+1)}\le \|x\|\le 2^{-j}\}.
$$
Since $\Delta_j\subset B(0,2^{-j})$ it can be covered using
$$ \cover_{X}(2^{-j},2^{-(j+3)})\le
    8^sM(\log 2)^2(1+|j|)^\beta(4+|j|)^\alpha
    \le \lemC(1+|j|)^{\alpha+\beta}:=M_j
$$
balls of radius $2^{-(j+3)}$, where $\lemC$ is independent of $j$.
We choose the centres $\{u^{(j)}_i\}_{i=1}^{M_j}$ of these balls
so that $ \|u^{(j)}_i\|\ge 2^{-(j+2)}.  $

Since $X$ is compact, $X\subset B(0,2^k)$ for some $k$
sufficiently large, and so
$$
\bigcup_{j=-k}^n\Delta_j= \{x\in X:\ \|x\|\ge 2^{-n}\}.
$$
Let $P_n$ be the orthogonal projection onto the linear subspace
$U_n$ spanned by the collection $\{u^{(j)}_i : j=-k,\ldots,n\hbox{
and }i=1,\ldots, M_j\}$.  Then the dimension of $U_n$ is bounded
by \mkconst{\lemD} $\lemD (1+n)^{\alpha+\beta+1}$ using the same
estimate as in \eqref{Cst}.  Moreover, for every $x\in \Delta_j$
there exists $u_i^{(j)}$ such that $x=u_i^{(j)}+v$ where $\|v\|\le
2^{-(j+3)}$. Since $\|P_n\|= 1$ and $\|P_nu\|=\|u\|$ for $u\in
U_n$, then
$$
\|P_nx\|=\|P_n(u_i^{(j)}+v)\|\ge\|P_nu_i^{(j)}\|-\|P_nv\|
    \ge 2^{-(j+2)}-2^{-(j+3)}
\ge \frac{1}{8} \|x\|.
$$

\end{proof}

\begin{blank}{Define
$$B_j=\{\, \omega\in \Re^{j} : \|\omega\|\le 1\,\}
$$
and let $\lambda_j$ be the uniform probability measure on $B_j$.
By Proposition \ref{mdev}, there are subspaces $U_k$ and
8-Lipschitz maps $\phi_k:U_k\to U_k^{\perp}$ such that
$\dist(X,G_{U_k}[\phi_k])\le 2^{-k}$ and $d_k=\dim U_k \le C
(1+k)^{\alpha+\beta+1}$. Since $U_k\sim \Re^{d_k}$ we may extend
the functionals $\omega\in B_{d_k}$ to $H$ by defining
$\omega(x)=\omega(P_k x)$ for $x\in H$. Define the probe set
\be{probeset}
Q=\{(l_1,\ldots,l_N):\ l_n=
    \frac{6}{\pi^2}\sum_{k=1}^\infty k^{-2}\phi_{nk}\
\mbox{with}\ \phi_{nk}\in B_{d_k}\}
\ee
Note that $Q$ is a compact subset of $S$;
moreover, any $L\in Q$ has Lipschitz
constant at most $1$.
Let $\mu$ be the probability measure on $Q$ that results from
choosing each $\phi_{nk}\in B_{d_k}$ 
independently with respect to the probability measure
$\lambda_{d_k}$. Before proving our main theorem we will provide a
few estimates concerning $\lambda_j$ and $\mu$.}
\end{blank}

Applying this lemma to $X-X$ there are subspaces $U_k$ with $\dim
U_k\le d_k:=c(1+k)^{\alpha+\beta+1}$ such that
$\|P_kz\|\ge\|z\|/8$ for all $z\in X-X$ with $\|z\|\ge 2^{-k}$.
Let $S_k$ denote the closed unit ball in $U_k$; clearly any
$\phi\in S_k$ induces a linear functional $L_\phi$ on $H$ via the
definition $L_\phi(u)=(\phi,u)$, where $(\cdot,\cdot)$ is the
inner product in $H$. Let $\zeta>0$ be fixed and define
$C_\zeta=1/\sum_{k=1}^\infty k^{-1-\zeta}$.
We now define the probe set
\be{probeset}
Q=\{(l_1,\ldots,l_N):\ l_n=L_{\phi_n}\ \mbox{where}\
\phi_n=C_\zeta\sum_{k=1}^\infty k^{-1-\zeta}\phi_{nk}\
\mbox{with}\ \phi_{nk}\in S_k\}.
\ee
We can identify $S_j$ with the unit ball $B_{d_j}$ in $\Re^{d_j}$,
and we denote by $\lambda_j$ the probability measure on $S_j$ that
corresponds to the uniform probability measure on $B_{d_j}$. We let
$\mu$ be the probability measure on $Q$ that results from choosing
each $\phi_{nk}$ randomly with respect to $\lambda_{d_k}$. Note that
$Q$ is a compact subset of $V$, and that all elements of $Q$ have
Lipschitz constant at most $\sqrt N$.

Before proving our main theorem we will prove a key estimate on
$\mu$. Although the argument is essentially the same as that in Hunt
\& Kaloshin (1999) our version is a little more explicit and we
include it here for completeness. The estimate relies on the
following simple inequality.

\begin{lemma}\label{mainlemA} If $x\in \Re^{j}$ and $\eta\in\Re$ then
$$
    \lambda_j\{\, \omega\in B_j :\ |\eta+(\omega\cdot x)|<\epsilon\,\}
        \le c j^{1/2} \epsilon |x|^{-1}.
$$
where $c$ is a constant that does not depend on $\eta$ or $j$.
\end{lemma}

\begin{proof}
Let $\hat x=x/|x|$. This follows immediately from estimate
\bea
    \lambda_j\{\, \omega\in B_j :\ |\eta+(\omega\cdot x)|<\epsilon\,\}
    &\le& \lambda_j\{\, \omega\in B_j
        : |\omega\cdot \hat x|<\epsilon |x|^{-1}\,\}\cr
    &=&\frac{\Omega_{j-1}}{\Omega_j}\,
        2\int_{0}^{\min(\epsilon |x|^{-1},1)}
        (1-\xi^2)^{(j-1)/2}\, \d\xi\cr
    &\le& \frac{\Omega_{j-1}}{\Omega_j}\, 2\epsilon\,|x|^{-1}
\eea
where $\Omega_j=\pi^{j/2}\Gamma(j/2+1)$ is the volume of the unit
ball in $\Re^j$. \end{proof}

\begin{lemma}\label{mainlemB} If $x\in H$ and $f\in V$ then
$$
    \mu\{\, L\in Q :\ |(L-f)(x)|<\epsilon \,\}
    \le c (d_k^{1/2}k^{1+\zeta}\epsilon \|P_k x\|^{-1})^N
$$
for every $k\in\N$ where $c$ is a constant independent of $f$ and
$k$.\end{lemma}

\begin{proof}
\def\J{{\cal J}}
Given $k\in \N$, let $\J$ be the index set $\J=\N\setminus \{ k \}$
and define $$B=\Big(\bigoplus_{j\in\J} B_{d_j}\Big)^N.$$
Given
$\alpha=((\alpha_{nj})_{j\in \J})_{n=1}^N\in B$ fixed, define
$$
    A_\alpha
    = \{\, (\phi_{nk})_{n=1}^N :\
        | (\eta_n+k^{-1-\zeta}\phi_{nk})(x)|<\epsilon
        \hbox{ for all } n\,\}
$$
where
$$\eta_n(x)=C_\zeta\sum_{j\in\J}
    j^{-1-\zeta}\alpha_{nj}(x) -f_n(x).$$
By Lemma \ref{mainlemA} there is a constant $c$ independent of
$\alpha$, $v$ and $k$ such that
$$ \lambda_{d_k}^N(A_\alpha)\le
        c (d_k^{1/2} k^{1+\zeta} \epsilon \|P_k x\|^{-1})^N.  $$
Let $P=\mu\{\, L\in Q :\ |(L-v)(x)|<\epsilon\,\}$. Then
$$
P\le\mu\{\, L\in Q:\ |(l_n-f_n)(x)|<\epsilon\ \hbox{for all}\ n\}.
$$
Let
$$
\Phi_N=
\{\,
        \big( (\phi_{nk})_{k=1}^\infty \big)_{n=1}^N
        :\ C_\zeta
    \Big|\sum_{k=1}^\infty k^{-1-\zeta} (\phi_{nk}-f_n)(x)\Big|
            <\epsilon,
        \forall n=1,\ldots, N
    \,\}
$$
Then by Fubini's theorem
\bea
    P&\le&\Big(\bigotimes_{j=1}^\infty \lambda_{d_j}\Big)^N
            \Phi_N\cr
    &=&\int_{\alpha\in B}\int_{\phi\in A_\alpha}
        \,\d \lambda_{d_k}^N(\phi)
        \, \d \Big(\bigotimes_{j\in\J} \lambda_{d_j}\Big)^N (\alpha) \cr
    &\le& \int_{\alpha\in B}
         c (d_k^{1/2} k^{1+\zeta} \epsilon \|P_k x\|^{-1})^N
        \, \d \Big(\bigotimes_{j\in\J} \lambda_{d_j}\Big)^N
        (\alpha)
    \cr
    &=&
        c (d_k^{1/2} k^{1+\zeta} \epsilon \|P_k x\|^{-1})^N.
\eea
This finishes the proof. \end{proof}

\begin{blank}{
\begin{lemma}\label{mainlemA}
If $x\in \Re^{j}$ and $\eta\in\Re$ then
$$
    \lambda_j\{\, \omega\in B_j :\ |\eta+\omega\cdot x|<\epsilon\,\}
        \le c j^{1/2} \epsilon \|x\|^{-1}.
$$
\end{lemma}
\begin{proof}
This follows immediately from estimate
$$
    \lambda_j\{\, \omega\in B_j :\ |\eta+\omega\cdot x|<\epsilon\,\}
    =\frac{\Omega_{j-1}}{\Omega_j}
        \int_{\max(\eta-\epsilon,-1)}^{\min(\eta+\epsilon,1)}
        (1-\xi^2)^{(j-1)/2}\, \d\xi
        \le \frac{\Omega_{j-1}}{\Omega_j}\, 2\epsilon
$$
where $\Omega_j=\pi^{j/2}/\Gamma(j/2+1)$ is the volume of the
unit ball in $\Re^j$.
\end{proof}

\begin{lemma}\label{mainlemB}
If $z\in X-X$ and $f\in V$ then
$$
    \mu\{\, L\in Q :\ |(L-f)(z)|<\epsilon \,\}
    \le c (d_k^{1/2}k^2\epsilon \|P_k z\|^{-1})^N
$$
for every $k\in\N$.
\end{lemma}

\begin{proof}
Let $P=\mu\{\, L\in Q :\ |(L-f)(z)|<\epsilon \,\}$. Then
\bea
\mu\{L\in Q:\ |Lx-f(z)|<\epsilon\}&\le&
\prod_{n=1}^N\,\nu\{l\in{\mathscr Q}:\ |l(x)-f_n(z)|<\epsilon\}
\eea
Since
$$
\nu\{l\in{\mathscr Q}:\ |l(z)-g(z)|<\epsilon\}\le\nu\{l\in{\mathscr
Q}:\ |l(z)|<\epsilon\}
$$
for any $g\in H^*$, it follows that
$$
P\le\big(\nu\{l\in{\mathscr Q}:\ |l(z)|<\epsilon\}\big)^N.
$$

Now,
\bea
\nu\{l\in{\mathscr Q}:\ |l(z)|<\epsilon\}
&=&\(\bigotimes_{k=1}^\infty\lambda_{d_k}\)\{\{\phi_k\}_{k=1}^\infty:\
\frac{6}{\pi^2}\left|\sum_{k=1}^\infty
k^{-2}(\phi_k,z)\right|<\epsilon\}\\
\eea

Now SOME WORK HERE let $k$ be such that $|z-P_kz|<\epsilon$. Then
since $l$ has Lipschitz constant no greater than one,
\bea
\nu\{l\in{\mathscr Q}:\ |l(z)|<\epsilon\}&\le&\nu\{l\in{\mathscr
Q}:\ |l(P_kz)|<2\epsilon\}\\
&\le&\(\bigotimes_{k=1}^\infty\lambda_{d_k}\)\{\{\phi_k\}_{k=1}^\infty:\
\frac{6}{\pi^2}\left|\sum_{k=1}^\infty
k^{-2}(\phi_k,z)\right|<\epsilon\}
\eea

Now, when $2^{-k}\le|z|\le 2^{-(k-1)}$  given $\epsilon=2^{-k}$ we
have

\bea
    P
    &\le & \Big(\bigotimes_{k=1}^\infty \lambda_{d_k}\Big)^N
    \{\,
        \big( (\phi_{nk})_{k=1}^\infty \big)_{n=1}^N
        :\ \frac{6}{\pi^2}
    \Big|\sum_{k=1}^\infty k^{-2}\, (\phi_{nk}-f_n,x)\Big|
            <\epsilon
        \mbox{ for all }n
    \,\}\\
    &\le & \Big(\bigotimes_{k=1}^\infty \lambda_{d_k}\Big)^N
    \{\,
        \big( (\phi_{nk})_{k=1}^\infty \big)_{n=1}^N
        :\ \frac{6}{\pi^2}
    \Big|\sum_{k=1}^\infty k^{-2}\, (\phi_{nk},x)\Big|
            <\epsilon
        \mbox{ for all }n
    \,\}
\eea

    &\le& \Big(\bigotimes_{k=1}^\infty \lambda_{d_k}\Big)^N
    \{\,
        \big( (\phi_{nk})_{k=1}^\infty \big)_{n=1}^N
        :\ |x\cdot(\phi_{nk}-f_n)|<{\pi^2}k^2\epsilon/6
        \mbox{ for all }n,k
    \,\}\\
    &=& \prod_{n=1}^N \prod_{k=1}^\infty \lambda_{d_k}
    \{\,
        \omega_{k}\in B_{d_k}
        :\ |x\cdot(\phi_{nk}-f_n)|<\pi^2k^2\epsilon/6
    \,\}.
\eea

Since $\lambda_k$ is a probability measure any of the terms in the
product can be dropped to obtain an upper bound.  Thus, for any
$k\in \N$ we obtain
\bea
    P&\le& \prod_{n=1}^N \lambda_{d_k}
    \{\,
        \omega_{k}\in B_{d_k}
        :\ |x\cdot(\phi_{nk}-f_n)|<\pi^2k^2\epsilon/6
    \,\}\\
    &\le&  c (d_k^{1/2}k^2\epsilon \|P_k x\|^{-1})^N
\eea
where the last inequality follows from Lemma \ref{mainlemA}.
\end{proof}

}\end{blank}

\subsection{Almost bi-Lipschitz embeddings}

We are now in a position to state and prove our main theorem, that
a compact subset $X$ of a Hilbert space with $X-X$ almost
homogeneous admits almost bi-Lipschitz linear embeddings into
finite-dimensional spaces. Unfortunately homogeneity of $X$ is not
automatically inherited by $X-X$: Olson (2002) exhibits an example
of a set $X$ with $d_{\rm A}(X)=0$ but $d_{\rm A}(X-X)=+\infty$
(for more see Section \ref{badones}).

\begin{theorem}\label{main}
Let $X$ be a compact subset of a Hilbert space $H$ such that $X-X$
is $(\alpha,\beta)$-almost homogeneous with $d_{\rm
A}^{\alpha,\beta}(X-X)<s<N$.  If
$$
\gamma>\frac{2+N(3+\alpha+\beta)+2(\alpha+\beta)}{2(N-s)}
$$
then a prevalent set of linear maps $f:H\rightarrow\Re^N$ are
injective on $X$ and, in particular, $\gamma$-almost bi-Lipschitz.
\end{theorem}

\ccount=0
\begin{proof}
First choose $\zeta>0$ in the definition of $Q$ small enough such that
\be{zetachoose}
\gamma>\frac{2+N(3+2\zeta+\alpha+\beta)+2(\alpha+\beta)}{2(N-s)}.
\ee

Since $\tau(X)\le\dF(X)\le\dF(X-X)\le\dA^{\alpha,\beta}(X-X)$ we can
apply Hunt \& Kaloshin's result (Theorem \ref{HK}, above) with
$\theta$ chosen so that
$$
    \theta> \frac {N(1+s/2)}{N-s}.
$$
to obtain a prevalent set $S_0$ of linear functions $f:H\to \Re^N$
\mkconst{\mainca} such that $f\in S_0$ implies there exists a
$\theta<1$ and $\mainca>0$ such that
\be{Szero}
    |f(x)-f(y)|\ge \mainca \|x-y\|^{\theta}
\qquad\hbox{for all}\qquad x,y\in X.
\ee
(We note here that the compactly supported probability measure used
in the definition of prevalence for $S_0$ differs from the measure
$\mu$ constructed in Section \ref{measure}, but is defined on the
same normed linear space $V$ of linear maps from $H$ to $\Re^N$). We
use this result to bootstrap a refined argument that makes use of
the stronger hypothesis that $\dA^{\alpha,\beta}(X-X)<\infty$.

Let $S_1$ be the subset of
$V$ consisting of those linear
functions $f:H\to \Re^N$ such that
$f\in S_1$ implies there exists $\delta>0$ such that
\be{Sone}
    |f(x)-f(y)|\ge \frac{\|x-y\|}{\slog(\|x-y\|)^\gamma}
\qquad\hbox{for all}\qquad \|x-y\|<\delta.
\ee
We now show that the set $S_1$ is also prevalent. Given $f\in V$,
let $K$ be the Lipschitz constant of $f$. We wish to show that
$\mu(f+S_1)=1$. This is equivalent to showing that $\mu(Q\setminus
(f+S_1))=0$.

Define the layers of $X-X$ by
\be{Zlayers}
Z_j=\{z\in X-X:\ 2^{-(j+1)}\le\|z\|\le 2^{-j}\}
\ee
and the set $Q_j$ of linear maps that fail to satisfy the required
continuity property\footnote{Strictly speaking the union of the
$Q_j$ form a set strictly larger than the complement of $S_1$.} for
some $z\in Z_j$ by
$$
Q_{j}=\{\,L\in Q:\ |(L-f)(z)|
        \le \Psi_{-\gamma}(2^{-j})
    \ \mbox{ for some }z\in Z_j\,\},
$$
where
$$
    \Psi_{-\gamma}(2^{-j})
    :=\frac{2^{-j}}{\sigma^\gamma\,\slog(2^{-j})^\gamma}
$$
and $\sigma$ is the constant occurring in (p4) in Lemma
\ref{SlogLem}. We now bound $\mu(Q_{j})$.

\mkconst{\maincb}
By assumption $\dA^{\alpha,\beta}(X-X)<s$, and so $Z_j$
can be covered by
\be{Mjest}
    M_j\le M\,
        \slog(2^{-j})^{\gamma
        s}\,\slog(2^{-j})^\beta\,\slog(\Psi_{-\gamma}(2^{-j}))^\alpha
    \le
    \maincb 
    (1+j)^{\alpha+\beta+\gamma s}
\ee
balls of radius $\Psi_{-\gamma}(2^{-j})$. Let the centres of these
balls be $z_i^{(j)}\in Z_j$ where $i=1,\ldots,M_j$. Given any
$z\in Z_j$ there is $z^{(j)}_i$ such that $\|z-z^{(j)}_i\|\le
\Psi_{-\gamma}(2^{-j})$. Thus
\bea
    |(L-f)(z)|&\ge& |(L-f)(z^{(j)}_i)|
            -|(L-f)(z-z^{(j)}_i)|\\
        &\ge& |(L-f)(z^{(j)}_i)|
            -(K+\sqrt N)\Psi_{-\gamma}(2^{-j})
\eea
implies
$$
    Q_{j}\subseteq
        \bigcup_{i=1}^{M_j}\{\,L\in Q:\ |(L-f)(z_i^{(j)})|
        \le (K+2\sqrt N)\Psi_{-\gamma}(2^{-j})\,\}.
$$

It follows, setting $k=j$ in Lemma \ref{mainlemB}, that
\bea\label{mainqest}
    \mu(Q_{j})&\le& \sum_{i=1}^{M_j}
        \mu \{\,L\in Q:\ |(L-f)(z_i^{(j)})|
        \le (K+2\sqrt N)\Psi_{-\gamma}(2^{-j})\,\}\\
        &\le& M_j \big(d_j^{1/2}j^{1+\zeta} (K+2\sqrt N)\Psi_{-\gamma}(2^{-j})
            \|P_j(z_i^{(j)})\|^{-1}\big)^N.
\eea
Now \eqref{Mjest} and Lemma \ref{stlem} imply that
$$
    \mu(Q_{j}) \le
        \maincb (1+j)^{\alpha+\beta+\gamma s}
        \big(d_j^{1/2}j^{1+\zeta}
        (K+2\sqrt N)2^{j+3}\Psi_{-\gamma}(2^{-j})\big)^N.
$$
\mkconst{\mainC} In particular (recall that $d_j\le C
(1+j)^{\alpha+\beta+1}$) there is a constant $\mainC>0$
independent of $j$ such that
$$
\mu(Q_{j})\sim \mainC
    j^{\alpha+\beta+\gamma s+N(\alpha+\beta+3+2\zeta-2\gamma)/2}
\quad \hbox{as} \quad j\to\infty.
$$
Since (\ref{zetachoose}) implies
$N(2\gamma-3-2\zeta-(\alpha+\beta))/2>1+\alpha+\beta+\gamma
s$, we have \mkconst{\mainD}
$$\sum_{j=1}^\infty\mu(Q_{j})<\mainD.$$

It follows from the Borel-Cantelli Lemma that $\mu$-almost every $L$
is contained in only a finite number of the $Q_j$; i.e.~there exists
a $J$ such that for all $j\ge J$, $ 2^{-(j+1)}\le\|z\|\le2^{-j}$
implies that $ |(L-f)(z)|\ge \Psi_{-\gamma}(2^{-j})$. It follows
from (p4) in Lemma \ref{SlogLem} that
$$
|(L-f)(z)| \ge
\sigma^\gamma\Psi_{-\gamma}(\|z\|)=\frac{\|z\|}{\slog(\|z\|)^\gamma}
\quad\hbox{for every}\quad \|z\|\le 2^{-J}.
$$
Thus $L-f\in S_1$ and so $L\in S_1+f$ for $\mu$-almost every $L$.

Define $S=S_0\cap S_1$.  Since the intersection of prevalent sets is
prevalent (Fact 3$'$ in Hunt et al. (1992)) $S$ is prevalent. Let
$f\in S$. Then there is $\mainca$ and $\delta$ such that both
\eqref{Szero} and \eqref{Sone} hold. Thus \mkconst{\mainE}
$$
    |f(x)-f(y)|\ge
    \mainE\,\frac{\|x-y\|}{\slog(\|x-y\|)^\gamma}
\qquad\hbox{for all}\qquad
    x,y\in X
$$
where $\mainE=\min\{1,\mainca \delta/\Psi_{-\gamma}(R)\}$ and $R>0$
is such that $X-X\subseteq B(0,R)$.
\end{proof}

Note that for a space $X$ with $X-X$ homogeneous,
i.e.~$\alpha=\beta=0$ in the above theorem, for any $\gamma>3/2$
we can choose $N$ large enough to obtain a $\gamma$-almost
bi-Lipschitz embedding into $\Re^N$.

We will prove a Banach space version of Theorem \ref{main} in
Section \ref{banachsec}. However, we delay this while, in the next
section, we consider in more detail almost homogeneity in a Hilbert
space.

\section{Lipschitz approximating dimension of Hilbert subsets and H\"older-Lipschitz embeddings}\label{HLE}

\ccount=0

The strong result of the previous section requires that $X-X$ is
almost homogeneous, while for a general almost homogeneous metric
space $(X,d)$ the embedding result of Theorem \ref{agone} only
provides a subset $f(X)$ of a Hilbert space that is itself almost
homogeneous.

Here we investigate further some of the properties of $f(X)$, and
are lead to define the `Lipschitz approximating dimension' and the
`Lipschitz deviation'. In particular we show that it is possible
to replace Hunt \& Kaloshin's thickness exponent with the
Lipschitz deviation.

\subsection{Further properties of the image
$f(X)$}\label{fXdelta0}

First we consider the almost bi-Lipschitz image $f(X)$ of a
compact almost homogeneous metric space $(X,d)$ in a Hilbert
space, as provided by Theorem \ref{agone}. We show that $f(X)$ can
be very well approximated by linear subspaces: it has `better than
zero' thickness.

As remarked after the proof of Theorem \ref{agone}, when $(X,d)$ is
compact the function $f$ defined by the simplified series
$$
f(x)=\sum_{j=-k}^\infty \frac{2^{-j}}{(1+|j|)^\delta M_j}\,
\phi_j(x)\otimes e_j
$$
still provides a $\gamma$-almost bi-Lipschitz embedding of $X$ into
a Hilbert space (choosing a $k$ such that $d(x_1,x_2)\le 2^{k}$ for
all $x_1,x_2\in X$). Now, for $n\in\N$ any element of $f(X)$ can be
approximated to within
$$
B \sum_{j=n+1}^{\infty}
        \frac{2^{-j}}{(1+|j|)^\delta}
   \ \le\ B\sum_{j=n+1}^{\infty} 2^{-j}
   \ \le\ B 2^{-n}
$$
\mkconst\Csttwo by an element of the subspace
$$
U=\bigoplus_{j=-k}^n \Re^{M_j}\otimes e_j,
$$
which has dimension \mkconst\Cstthr
\be{Cst}\sum_{j=-k}^n M_j\le (n+k+1)C
        (1+n)^{\alpha+\beta}
    \le \Csttwo (1+n)^{\alpha+\beta+1}.
\ee
Here $\Csttwo$ depends on $C,$ $k$ and the constants in Lemma
\ref{SlogLem} but is independent of $n$. It follows that
\be{delta0}
  d(f(X),\epsilon)\le \Cstthr
    \big[\log(\e+1/\epsilon)\big]^{\alpha+\beta+1}.
\ee
One consequence of this inequality is that the thickness exponent
of $f(X)$ is zero, but (\ref{delta0}) is significantly stronger
than this.

\subsection{The Lipschitz deviation}

Inspired by the quantity $d(X,\epsilon)$ used to define the
thickness we now introduce a more general quantity, the
$m$-Lipschitz deviation: we denote by $\delta_{m}(X,\epsilon)$ the
smallest dimension of a linear subspace $U$ such that
$$
\dist(X,G_U[\phi])<\epsilon
$$
for some $m$-Lipschitz function $\phi:U\rightarrow U^\perp$,
$$
\|\phi(u)-\phi(v)\|\le m\|u-v\|\qqfa u,v\in U,
$$
where $U^\perp$ is the orthogonal complement of $U$ in $H$. We
will write $G_U[\phi]$ for the graph of $\phi$ over $U$:
$$
G_U[\phi]=\{\,u+\phi(u):\ u\in U\,\}.
$$
Clearly $\delta_0(X,\epsilon)=d(X,\epsilon)$.

In Section \ref{fXdelta0} we showed that for the almost
bi-Lipschitz embedding $f(X)$ of an almost homogeneous metric
space into a Hilbert space
$$
\delta_0(f(X),\epsilon)\le\Cstthr
    \big[\log(\e+1/\epsilon)\big]^{\alpha+\beta+1}.
$$
We now show that Lemma \ref{stlem} implies a bound of a similar
form on $\delta_8(X,\epsilon)$ for any subset of a Hilbert space
with $X-X$ almost homogeneous.

\begin{proposition}\label{mdev}
Let $X$ be a compact subset with the set of differences $X-X$
$(\alpha,\beta)$-almost $(M,s)$-homogeneous. Then there exists a
sequence of linear subspaces $U_k$ with $\dim U_k\le C
(1+k)^{\alpha+\beta+1}$ and $U_{k+1}\supseteq U_k$, and
$8$-Lipschitz functions $\phi_k:U_k\rightarrow U_k^\perp$ such
that
$$
\dist(X,G_{U_k}[\phi_k])\le 2^{-k}.
$$
In particular
$$
\delta_8(X,\epsilon)\le
K\big[\log(\e+1/\epsilon)\big]^{\alpha+\beta+1}.
$$
\begin{blank}{
$\mdev_8(X)=1$. Furthermore $X$ is contained in the graph of an
$8$-Lipschitz function over $U$. }\end{blank}
\end{proposition}

\begin{proof}
Applying Lemma \ref{stlem} to $X-X$ we obtain a nested sequence of
linear subspaces for which
$$
\frac{1}{8}\|x-y\|\le\|P_kx-P_ky\|\le \|x-y\|\qqfa x,y\in X\qwith
\|x-y\|\ge 2^{-k},
$$
where $P_k$ is the orthogonal projection onto $U_k$.

Define $\phi_k:U_k\rightarrow U_k^\perp$ as follows. Let $N_k$ be
a maximal $2^{-k}$ net in $(X,d)$ and set $ \phi_k(P_kx)=
\Qk x$ for $x\in N_k$.
Given $P_k x,P_ky\in P_k N_k$ we have
$$
    \|\phi_k(P_kx)-\phi_k(P_ky)\|
        \le \|\Qk (x-y)\|\le \|x-y\|
        \le 8 \|P_kx -P_k y\|.
$$
Therefore $\phi_k: P_k N_k\to U_k^\perp$ is a $8$-Lipschitz
function. Now, extend this $\phi_k$ to a $8$-Lipschitz function
$U_k\rightarrow U_k^\perp$.

Since $N_k\subset G_{U_k}[\phi_k]$ then any point of $X$ lies
within $2^{-k}$ of $G_{U_k}[\phi_k]$. Thus
$$
\delta_8(X,2^{-k})\le \lemC (1+k)^{\alpha+\beta+1}
$$
and the result follows.\end{proof}

We now show that this argument can be reversed, i.e.~that the
results of Lemma \ref{stlem} and Proposition \ref{mdev} are
essentially equivalent.

\begin{proposition}\label{reverse}
Suppose that $X$ is a compact subset of a Hilbert space $X$. For
any $m\ge0$ let $\{U_k\}_{k=1}^\infty$ be a sequence of linear
subspaces such that for each $U_k$ there exists an $m$-Lipschitz function
$\phi_k:U_k\rightarrow U_k^\perp$ with
$$
\dist(X,G_{U_k}[\phi_k])\le 2^{-k}.
$$
Then there exists an integer $n$ and a constant $c_m>0$ (which
depends on $m$ but is independent of $k$) such that for every $k$
$$
\|P_{k+n}(x_1-x_2)\|\ge c_m\,\|x_1-x_2\|\qfa x,y\in X\qwith
\|x_1-x_2\|\ge 2^{-k}.
$$
\end{proposition}

\begin{proof}
First note that for any $x\in H$ we have
$$
\dist(x,G_{U_k}[\phi_k])^2=\inf_{u\in
U_k}\bigl(\|P_kx-u\|^2+\|\Qk x-\phi_k(u)\|^2\bigr)
$$
and since for any $u\in U_k$ we have
\bea
\|\Qk x-\phi_k(P_kx)\|^2&=&\|\Qk x-\phi_k(u)+\phi_k(u)-\phi_k(P_kx)\|^2\\
&\le& 2\|\Qk x-\phi_k(u)\|^2+2\|\phi_k(u)-\phi_k(P_kx)\|^2\\
&\le& 2\|\Qk x-\phi_k(u)\|^2+2m^2\|u-P_kx\|^2\\
&\le& l_m^2\bigl(\|P_kx-u\|^2+\|\Qk x-\phi_k(u)\|^2\bigr),
\eea
where $l_m^2=2\max(1,m^2)$, it follows that for $x\in X$
\be{G2d}
\|\Qk x-\phi_k(P_kx)\|\le l_m\,\dist(x,G_{U_k}[\phi_k])\le l_m 2^{-k}
\ee

Now suppose that $x_1,x_2\in X$ with
$$
\|x_1-x_2\|\ge 2^{-k}.
$$
Let $n$ be the smallest integer
such that $3 l_m\le 2^n$ and
set
$$
\tilde x_j=P_{k+n}x_j+\phi_{k+n}(P_{k+n}x_j)\qquad\hbox{for}\qquad j=1,2.
$$
Clearly, $P_{k+n}(x_1-x_2)=P_{k+n}(\tilde x_1-\tilde
x_2)$.
Furthermore, it follows from (\ref{G2d}) that
$|x_j-\tilde x_j|\le 2^{-k}/3$ for $j=1,2$.
Therefore, $|\tilde x_1-\tilde x_2|\ge|x_1-x_2|/3$.

Now, since $\tilde x_1,\tilde x_2\in G_{U_{k+n}}[\phi_{k+n}]$,
\bea
\|P_{k+n}\tilde x_1-P_{k+n}\tilde x_2\|^2&=&\|\tilde x_1-\tilde
x_2\|^2
-\|\phi_{k+n}(P_{k+n}\tilde x_1)-\phi_{k+n}(P_{k+n}\tilde x_2)\|^2\\
&\ge&\|\tilde x_1-\tilde x_2\|^2-m^2\|P_{k+n}(\tilde x_1-\tilde
x_2)\|^2,
\eea
and so
$$
\|P_{k+n}(x_1-x_2)\|=\|P_{k+n}(\tilde x_1-\tilde
x_2)\|\ge\frac{\|\tilde x_1-\tilde
x_2\|}{\sqrt{1+m^2}}\ge\frac{\|x_1-x_2\|}{3\sqrt{1+m^2}}.
$$
\end{proof}

\subsection{Almost homogeneous subsets of a Hilbert space}

If we assume only the almost homogeneity of $X$, rather than of
$X-X$, we can apply a simplified variant of the argument of
Theorem \ref{main} to obtain the following minor improvement to
the embedding theorem of Hunt \& Kaloshin (under our stronger
hypothesis). For a zero thickness set $X$ with $\dF(X)\le d$ they
obtain an upper limit of $N/(N-2d)$ for the H\"older exponent,
while under the assumption that $\dA^{\alpha,\beta}(X)\le s$ we
obtain $(N-s)/(N-2s)$ as the upper limit. Note that we replace any
assumption on the thickness by (\ref{deltamassump}), which in
particular is satisfied by the almost bi-Lipschitz embedding
$f(X)$ of an almost homogeneous metric space with $m=0$ (see
\ref{delta0}).

\begin{theorem}\label{subT}
Suppose that $X$ is a compact subset of a Hilbert space $H$ with
$d_{\rm A}^{\alpha,\beta}(X)< s$ and that for some $m>0$,
$\sigma\ge0$,
\be{deltamassump}
\delta_m(X,\epsilon)\le K[\log(\e+1/\epsilon)]^\sigma.
\ee
Then for any integer $N>2s$, if $\theta>(N-s)/(N-2s)$ there is a
prevalent set $S$ of linear maps $f:H\rightarrow\Re^N$ such that for
every $f\in S$ there exists $c>0$ such that
\be{newHolder}
|f(x)-f(y)|\ge c\|x-y\|^\theta
\quad\mbox{for all}\quad x,y\in X.
\ee
\end{theorem}

\ccount=0
\begin{proof}
Set
$$d_j=\delta_m(X,2^{-j})\le K\big[\log(\e+2^j)\big]^\sigma$$
and define $Q$ as in \eqref{probeset} with
$\zeta=1$.
Define the layers $Z_j$ as in \eqref{Zlayers}
and
$$
Q_j=\{\, L\in Q:\ |(L-v)(z)|\le 2^{-j\theta}\mbox{ for some }z\in Z_j\,\}.
$$
\mkconst{\subca}%
Let $R>0$ be chosen so large that $X\subset B(0,R)$. Cover $X$ by
\bea
\cover_X(R,2^{-(j+1)\theta}) &\le& M
\(\frac{R}{2^{-(j+1)\theta}}\)^s\slog(R)^\beta\slog(2^{-(j+1)\theta})^\alpha\\
    &\le& \subca 2^{j\theta s}(1+j\theta)^\alpha
\eea
balls of radius $2^{-(j+1)\theta}$ centred at points $x_i\in X$.
Denote these as
$$
 X_i=\{\, x\in X: \|x-x_i\|<2^{-(j+1)\theta}\,\}.
 $$
Now consider the larger balls
$$
B_i=\{\,y\in X:\ \|x_i-y\|\le 2^{-(j+1)\theta}+2^{-j}\,\}.
$$
Cover each of these balls by at most
\mkconst{\subcb}
\bea
   && \cover_X(2^{-(j+1)\theta}+2^{-j},2^{-(j+1)\theta})\\
        &&\qquad\le
        M\(1+2^{(j+1)\theta-j}\)^s\,\slog(2^{-(j+1)\theta}+2^{-j})^\beta\,\slog(2^{-(j+1)\theta})^\alpha\\
    &&\qquad\le \subcb 2^{j(\theta-1) s}(1+j)^\beta(1+j\theta)^\alpha
\eea
balls of radius $2^{-(j+1)\theta}$.
Since
$$
    Z_j = \bigcup_i \bigcup_{x\in X_i}
        \{\, x-y : 2^{-(j+1)}<\|x-y\|<2^{-j}\,\}
    \subseteq \bigcup_i (X_i-B_i)
$$
it follows that $Z_j$ can be covered by
$$
M_j=\subca\subcb 2^{js(2\theta-1)}(1+j\theta)^{2\alpha}(1+j)^\beta
$$
balls of radius $2^{-j\theta}$.  Let $z_i^{(j)}$ denote the centres
of these balls.

Applying similar estimates as in the proof of Theorem \ref{main}
(these rely on Proposition \ref{reverse} to ensure that
$\|P_kz_i^{(j)}\|\ge c\|z_i^{(j)}\|$ for some $c>0$) shows that
\mkconst{\subcc}
$$
    \mu(Q_{j})\sim
        2^{js(2\theta-1)}j^{2\alpha+\beta}[j^{2+\sigma}2^{j(1-\theta)}]^N
        \quad\hbox{as}\quad j\to\infty.
$$
Thus $\sum \mu(Q_j)$ converges provided that $\theta>(N-s)/(N-2s)$.
The argument is now concluded as in Theorem \ref{main}.
\end{proof}

By combining this with Proposition \ref{agone} we obtain the
following H\"older-Lipschitz embedding result for homogeneous
metric spaces (cf.~Lemma 9.1 in Foias and Olson (1996) which has a
similar result for spaces with finite upper box-counting
dimension).

\begin{corollary}
Let $(X,d)$ be an almost homogeneous metric space with $d_{\rm
A}^{\alpha,\beta}<s$. If $N>2s$ and $\theta>(N-s)/(N-2s)$ there
exists a map $\phi:(X,d)\rightarrow\Re^N$ such that
$$
c^{-1}\,d(x,y)^{\theta}\le|\phi(x)-\phi(y)|\le c\,d(x,y)\qqfa
x,y\in X.
$$
\end{corollary}

Of course one can prove finite-dimensional versions of Theorems
\ref{main} and \ref{subT} using very similar techniques.

\subsection{The Lipschitz deviation}

It is interesting that our argument shows that for any fixed $m>0$
the thickness exponent in the statement of Theorem \ref{HK} can be
replaced by the {\it $m$-Lipschitz deviation}, $\dev_m(X)$, which
we define by analogy with the thickness exponent
(cf.~(\ref{thickness}))
$$
\dev_m(X)=\limsup_{\epsilon\rightarrow0}\frac{\log\delta_m(X,\epsilon)}{-\log\epsilon}.
$$
We note that $\dev_m(X)\le\tau(X)$ and that this gives an
indication of why the thickness exponent can be expected to play a
r\^ole in determining the H\"older exponent in (\ref{newHolder}).
We state without proof:

\begin{theorem} Let $X$ be a compact
subset of a Hilbert space $H$, $D$ an integer with $D>\dF(X-X)$,
and let $\dev_m(X)$ be the $m$-Lipschitz deviation of $X$. If
$\theta$ is chosen with
$$\theta>\frac{D(1+\dev_m(X)/2)}{D-\dF(X-X)}$$
then for a prevalent set of linear maps $L:B\rightarrow\Re^D$
there exists a $c>0$ such that
$$
c\|x-y\|^\theta\le|Lx-Ly|\le \|L\|\|x-y\|
    \quad\mbox{for all}\quad x,y\in X;
$$
in particular these maps are injective on $X$.
\end{theorem}

\section{Embedding subsets $X$ of Banach spaces with $X-X$
homogeneous}\label{banachsec}

In this section we extend the Hilbert space result to cover subsets
of Banach spaces. In particular this enables us to prove a new
almost bi-Lipschitz embedding result for a class of metric spaces.

The key point is, of course, that enough of Lemma \ref{stlem} can be
salvaged to follow a very similar proof:

\begin{lemma}\label{banachone}
Let $X$ be an $(\alpha,\beta)$-almost $(M,s)$-homogeneous subset of
a Banach space $B$. Then there exists a nested sequence of subsets
$U_{n+1}\supseteq U_n$ such that
$$
\dim U_n\le C(1+n)^{\alpha+\beta+1}
$$
and
$$
\dist(x,U_n)\le \frac{1}{4}\|x\|\qqfa \|x\|\ge 2^{-n}
$$
\end{lemma}

In particular, if we apply this lemma to $Z=X-X$, there exists a
nested sequence of linear subspaces of $B$, $U_k\subseteq U_{k+1}$
such that given $z\in X-X$ with $\|z\|\ge 2^{-n}$ there exists a
point $\tilde z\in U_n$ such that
$$
\|z-\tilde z\|\le \frac{1}{4}\|z\|\qqand\|\tilde
z\|\ge\frac{3}{4}\|z\|.
$$

We now let $S_k$ denote the closed unit ball in the dual of $U_k$,
and denote by $S_k^E$ an isometric embedding of $S_k$ into $B^*$,
whose existence is guaranteed by the Hahn-Banach theorem. We then
define our probe set $Q$ as
$$
Q=\{(l_1,\ldots,l_N):\ l_n= C_\zeta \sum_{k=1}^\infty
k^{-1-\zeta}\phi_{nk}\ {\rm with}\ \phi_{nk}\in S_k^E\}.
$$
Choosing a basis for $S_k$ we identify $S_k$ with a convex set
$U_k\subset\Re^{d_k}$, and induce a probability measure on $S_k$
(and hence on $S_k^E$) via the uniform probability measure on $U_k$.

We now outline the proof of the following result:

\begin{theorem}\label{mainbanach}
Let $X$ be a compact subset of a Hilbert space $B$ such that $X-X$
is $(\alpha,\beta)$-almost homogeneous with $d_{\rm
A}^{\alpha,\beta}(X-X)<s<N$.  If
$$
\gamma>\frac{1+N(2+\alpha+\beta)+(\alpha+\beta)}{N-s}
$$
then a prevalent set of linear maps $f:B\rightarrow\Re^N$ are
injective on $X$ and, in particular, $\gamma$-almost bi-Lipschitz.
\end{theorem}

\begin{proof} The proof proceeds identically to that of Theorem \ref{main}
until we have to estimate
$$
\mu\(L:\ |(L-f)z_i^{(j)}|\le c\Psi_{-\gamma}(2^{-j})\).
$$
We can now follow the argument from Hunt \& Kaloshin (1999), with
some small changes---we only highlight these here. In our case we
know that there exists a point $\tilde z_i^{(j)}\in U_j$ such that
$$
\|\tilde z_i^{(j)}-z_i^{(j)}\|\le \frac{1}{4}\|z_i^{(j)}\|.
$$
It follows that there exists a $\psi\in S_j$ such that
$$
\psi(z_i^{(j)})\ge\|z_i^{(j)}\|-\|z_i^{(j)}-\tilde
z_i^{(j)}\|\ge\frac{3}{4}\|z_i^{(j)}\|\ge 3\cdot 2^{-(j+3)}.
$$
We can then follow Hunt \& Kaloshin's argument to show that
$$
\mu\(L:\ |(L-f)z_i^{(j)}|\le
c\Psi_{-\gamma}(2^{-j})\)\le\(j^{1+\zeta}d_j2^{j+3}\Psi_{-\gamma}(2^{-j})\)^N,
$$
and the proof is completed exactly as in the Hilbert space case,
noting that we now have a factor of $d_j$ rather than only
$d_j^{1/2}$.
\end{proof}

One significant consequence of extending the result to Banach spaces
is it allows for a new result for metric spaces via the Kuratowski
isometric embedding of $(X,d)$ into $L^\infty(X)$: choosing an
arbitrary point $x_0\in X$, this is given by
\be{Kuratowski}
 x\mapsto
\rho_x,\qquad{\rm where}\qquad\rho_x(y)=d(x,y)-d(x_0,y)
\ee
(see Heinonen, 2003, for example). In this way we can attach meaning
to $X-X$ for an arbitrary metric space $(X,d)$, i.e.~
\be{XminusX}
X-X=\{f\in L^\infty(X):\ f=d(x,\cdot)-d(y,\cdot),\ x,y\in X\}.
\ee
We then have the following result:

\begin{theorem}\label{metricone}
Let $(X,d)$ be a compact metric space such that $X-X$ is an almost
homogeneous subset of $L^\infty(X)$. Then there exists an injective
almost bi-Lipschitz map $f:(X,d)\rightarrow\Re^N$.
\end{theorem}

\begin{proof}
Denote by ${\mathscr F}:(X,d)\rightarrow L^\infty(X)$ the isometric
embedding in (\ref{Kuratowski}). Then ${\mathscr F}(X)$ is isometric
to $(X,d)$, while the set of differences ${\mathscr F}(X)-{\mathscr
F}(X)$ is almost homogeneous by assumption. The existence of an
injective almost bi-Lipschitz embedding of ${\mathscr F}(X)$ into
$\Re^N$, which follows from the Banach space version of our main
theorem, immediately implies the existence of the same type of
embedding for $(X,d)$ into $\Re^N$.
\end{proof}

\section{The relationship between $\dA^{\alpha,\beta}(X)$ and
$\dA^{\alpha,\beta}(X-X)$}\label{badones}

In this section we give some results relating the homogeneity of $X$
and $X-X$.
First, we give an example of a set $X$ for which
$d_{\rm A}(X)=0$ but $d_{\rm A}(X-X)=+\infty$. It is easy to show
that the set
\be{Xstar}
X^*=\{a_ne_n:\ a_n=4^{-(2^j)}, n=2^{j-1},\ldots,2^j-1\},
\ee
where $e_n$ is an orthonormal basis of a Hilbert space $H$, has
$d_{\rm A}(X^*)=+\infty$. Note that $|a_n|\le 4^{-n}$ for all $n$.
Consider now the subset $X$ of $H\times H$ defined by
$$
X=\{(4^{-n}e_n,a_ne_n)\}_{n=1}^\infty\cup\{(4^{-n}e_n,0)\}.
$$
A simple argument shows that $d_{\rm A}(X)=0$, while $X-X$
contains a copy of $X^*$ and so $d_{\rm A}(X-X)=\infty$.

This negative result appears to be in some ways typical
for almost homogeneous sets as well, as we
will now show. We begin with a preparatory lemma.

\begin{lemma}\label{notalg}
The orthogonal sequence with algebraic decay
$$X^*=\{\, b_n e_n : b_n \sim \epsilon n^{-\gamma}\,\}$$
where $\epsilon,\gamma>0$
has $\dA^{\alpha,\beta}(X^*)=+\infty$
for any $\alpha,\beta\ge 0$.
\end{lemma}

\begin{proof}
Let $n_0$ be chosen so large that
$$
    \epsilon (2n)^{-\gamma} < |b_n| < \epsilon(n/2)^{-\gamma}
 \qquad\hbox{for}\qquad n>n_0.
$$
Let $r_n=\epsilon(n/2)^{-\gamma}$ and $\rho_n=\epsilon
(4n)^{-\gamma}$. Suppose, for a contradiction, that
$\dA(X^*)<s<\infty$. Then there exists an $M\ge 1$ such that
\be{OSADA}
    \cover(r_n,\rho_n)
\le M\Big(\frac{r_n}{\rho_n}\Big)^s
\slog(r_n)^\beta\slog(\rho_n)^\alpha.
\ee
On the other hand,
$$
    B(0,r_n) \supseteq \{\, b_k e_k : n<k \le 2n \},
$$
where the points $b_ke_k$ with $n<k\le 2n$ are each a distance
greater than $|b_k|>\epsilon (4n)^{-\gamma}$ apart from each other.
Therefore,
\be{OSADB}
    \cover(r_n,\rho_n)
        \ge {\rm card} \big(\{\, b_k e_k : n<k\le 2n \}\big) = n.
\ee
Combining inequality (\ref{OSADA}) with (\ref{OSADB}) and applying
(p1) of Lemma 2.1 we obtain
$$
    n \le M 8^{\gamma s}
    \big(\log 2+|\log \epsilon(n/2)^{-\gamma}| \big)^\beta
    \big(\log 2+|\log \epsilon(4n)^{-\gamma}| \big)^\beta.
$$
Letting $n\to \infty$ yields a contradiction, and so
$\dA^{\alpha,\beta}(X^*)=\infty$.
\end{proof}

\begin{lemma}\label{feta}
Given two unit vectors $v,w\in H$ set $e_1=v$ and choose
$\alpha\in \Re$ and a unit vector $e_2$
such that $e_1\cos\alpha-e_2\sin\alpha=w$ and
$\cos\alpha=(v,w)$. Note that $e_2$ is orthogonal to $e_1$.
Extend $\{e_1,e_2\}$ to a basis for $H$, and define the rotation
$$
R_x=\left(\begin{array}{cc}\cos(\alpha\psi(x))&\sin(\alpha\psi(x))\\
-\sin(\alpha\psi(x))&\cos(\alpha\psi(x))\end{array}\right)\oplus{\rm
id},
$$
where $\psi:H\rightarrow\Re$ is a fixed $C^\infty$ function such
that
$$
\psi(x)=\{\begin{array}{cl}
    0&\hbox{if }\ \|x\|\le 3/4\hbox{ or }\|x\|\ge 2,\\
    1&\hbox{if }\ \|x\|=1.
\end{array}\right.
$$
Let $f(x)=R_x x$.
Then $f\in C^\infty$ and $f(v)=w$.
Moreover, $f_\eta(x)=\eta^{-1} f(\eta x)$ is uniformly
bi-Lipschitz continuous for $\eta>0$ and
different from the identity only for
$x\in H$ such that
$(3/4) \eta^{-1}<\|x\|<2\eta^{-1}$.
\end{lemma}

\begin{proof}
By construction $f\in C^\infty$, $f(v)=w$ and $f(x)=x$ for
$\|x\|\le 3/4$ or $\|x\|\ge 2$.  Rescaling shows that $f_\eta(x)$
is different from the identity only for
$(3/4)\eta^{-1}<\|x\|<2\eta^{-1}$. We now show that $f_\eta(x)$ is
uniformly bi-Lipschitz continuous for $\eta>0$.

Let $x,y\in H$ with $\|x\|\le \|y\|$. If $\|x\|\ge 2\eta^{-1}$ then
$f_\eta(x)=x$ and $f_\eta(y)=y$, so we consider only the case
$\|x\|<2\eta^{-1}$. Then
\bea
\|f_\eta(x)-f_\eta(y)\|&=&\|R_{\eta x}x-R_{\eta y}y\|\\
&\le&\|(R_{\eta x}-R_{\eta y})x\|+\|R_{\eta y}(x-y)\|\\
&\le&\|R_{\eta x}-R_{\eta y}\|\|x\|+\|R_{\eta y}\| \|x-y\|\\
&\le&2\eta^{-1}\|R_{\eta x}-R_{\eta y}\|+\|x-y\|.
\eea
Since
\par\noindent\vbox{
$$
\|R_{\eta x}-R_{\eta
y}\|=\left\|\left(\begin{array}{cc}
\cos(\alpha\psi(\eta x))-\cos(\alpha\psi(\eta y))
&\sin(\alpha\psi(\eta x))-\sin(\alpha\psi(\eta y))\\
-\sin(\alpha\psi(\eta x))+\sin(\alpha\psi(\eta y))
&\cos(\alpha\psi(\eta x))-\cos(\alpha\psi(\eta y))
\end{array}\right)\right\|
$$
$$
\le C_1\alpha\eta\|x-y\|:=C_2\eta\|x-y\|,\hfill
$$}
it follows that
$$
\|f_\eta(x)-f_\eta(y)\|\le (2C_2+1)\|x-y\|
$$
where the Lipschitz constant $2C_2+1$ does not depend on $\eta$.
Since $f_\eta$ is injective with inverse $f_\eta^{-1}$ formed by
the same construction but with the roles of $v$ and $w$ reversed
we obtain the same bound for $\|f^{-1}_\eta(x)-f^{-1}_\eta(y)\|$.
\end{proof}

\begin{proposition}\label{nasty}
Let $X$ be a connected subset of a Hilbert space $H$ that contains
more than one point. Then there exists a $C^\infty$ bi-Lipschitz
map $\phi:H\rightarrow H$ such that
$$
\dA^{\alpha,\beta}(\phi(X)-\phi(X))=+\infty
$$
for every $\alpha,\beta\ge 0$.
Furthermore $\phi$ may be chosen such that $\dist_{\rm
H}(\phi(X),X)$ is arbitrarily small.
\end{proposition}

\begin{proof}
Since $X$ contains more than one point,
there exist two disjoint
balls $B(x_1,R)$ and $B(x_2,R)$ of radius~$R>0.$
Moreover, since $X$ is connected, then there are
points $x_{2+i}\in X$ for $i=1,2$
such that $\|x_{2+i}-x_i\|=R/4$.
Thus, the four balls $B(x_i,R/8)$
with $x_i\in X$ for $i=1,\ldots, 4$ are disjoint.
Moreover
$$\bigcup_{i=1}^4 B(x_i,R/8) \subseteq
    \bigcup_{i=1}^2 B(x_i,3R/8).
$$
Recursively define nested families of disjoint balls such that
$$\bigcup_{i=1}^{2^{j+1}} B(x_i,R8^{-j}) \subseteq
    \bigcup_{i=1}^{2^{j}} B(x_i,3R8^{-j}).
$$

For $j=0,1,2,\ldots$ and $i=1,\ldots,2^{j+1}$
let $a_j=(1/2)R8^{-j}$ and $e_{ij}=e_{2^{j+1}-2+i}$
where $e_i$ is an orthonormal basis of $H$.
Choose the points $y_{ij}\in B(x_i,R8^{-j})$
such that $\|x_i-y_{ij}\|=a_j$.
Further define
$$
g_{ij}(x)=x_i+f_{\eta}(x-x_i),
$$
where $f_\eta$ is the function given in Lemma \ref{feta} for
$v=(y_{ij}-x_i)/a_j$, $w=e_{ij}$
and $\eta=1/a_j$.
If $\|x-x_i\|\ge 2a_j=R8^{-j}$ or
$\|x-x_i\|\le (3/4)a_j=3R8^{-j-1}$ then
$f_\eta(x-x_i)=x-x_i$ and $g_{ij}(x)=x$.
Therefore the function $g_{ij}$ is
$C^\infty$, bi-Lipschitz and different from the identity only
on the annulus $B(x_i,R8^{-j})\setminus B(x_i,3R8^{-j-1})$.
Moreover, by construction we have
$$
g_{ij}(y_{ij})=x_i+f_\eta(y_{ij}-x_i)=x_i+a_if(v)=x_i+ a_ie_{ij}.
$$
Set
$$\phi(x)=\sum_{j=0}^\infty\sum_{i=1}^{2^{j+1}} g_{ij}(x).$$
Since the
$g_{ij}$ are different from the identity only on disjoint sets
and the bi-Lipschitz
constant of $f_\eta$ is independent of $\eta$,
then the map $\phi$ is a
bi-Lipschitz $C^\infty$ map of $H$ onto $H$. Since
$\phi(X)-\phi(X)$ contains
\bea
    &&\{\, a_j e_{ij} :j=0,1,2,\ldots\hbox{ and }
            i= 1,\ldots,2^{j+1}\, \} \\
        &&\qquad= \{\, b_n e_n :
           b_n=(1/2)R8^{-j} ,\ n=2^{j+1}-1,\ldots, 2^{j+2}-2
     \,\}
\eea
where $4R/(n+2)^3\le b_n\le 4R/(n+1)^3$, then
 $b_n\sim 4R n^{-3}$ and hence Lemma~\ref{notalg}
implies $\dA^{\alpha,\beta}(\phi(X)-\phi(X))=\infty$.

Finally, note that $\dist_{\rm H}(\phi(X),X)$
may be made arbitrarily small by taking $R>0$
sufficiently small in step one.
\end{proof}

A consequence of this result is that it is not necessary for $X-X$
to be homogeneous in order to obtain a bi-Lipschitz embedding of
$X$ into some $\Re^k$. Indeed, any set $X$ that can be so embedded
has a bi-Lipschitz image that has
$\dA^{\alpha,\beta}(X-X)=\infty$. However, it may still be the
case that $X-X$ has to be homogeneous in order to obtain a {\it
linear} bi-Lipschitz embedding as in Theorem \ref{main}.

On a more positive note, if $X$ is an orthogonal sequence then
homogeneity of $X$ does imply homogeneity of $X-X$.

\begin{lemma}\label{XmX}
Let $X=\{\,x_j\}_{j=1}^\infty$ be an orthogonal sequence in $H$.
If $d_{\rm A}(X)<+\infty$ then $d_{\rm A}(X-X)\le 2d_{\rm A}(X)$.
\end{lemma}

\begin{proof}
Suppose that $X$ is $(M,s)$-homogeneous. We write
$B_X(r,x)=B(r,x)\cap  X$, and consider a ball
$B=B_{X-X}(r,x-y)\subseteq X-X$ of radius $r$ centred at $x-y\in
X-X$. Since $B\subseteq B_{X-X}(\rho,0)\cup \big(B\setminus
\{0\}\big)$, we need only cover $B\setminus \{0\}$.

Suppose that $x=y$, so that $B=B_{X-X}(r,0)$.  Let $a-b\in
B\setminus\{0\}$. Then $a\ne b$ and therefore $a$ is orthogonal to
$b$.  It follows that
$$
    \big\|(a-b)-(x-y)\big\|^2=\|a\|^2+\|b\|^2 < r^2.  $$
Hence $a,b\in B_X(r,0)$, and consequently
$$
    B\setminus\{0\} \subseteq B_X(r,0)-B_X(r,0).$$
Cover $B_X(r,0)$ with $M(2r/\rho)^s$ balls $B_X(\rho/2,a_i)$ of
radius $\rho/2$ centred at $a_i\in X$.  Then
\bea
    \bigcup_{i,j} B_{X-X}(\rho,a_i-a_j)
    &\supseteq& \bigcup_i B_X(\rho/2,a_i)-\bigcup_j B_X(\rho/2,a_j)\\
    &\supseteq& B_X(r,0)-B_X(r,0)\supseteq B_{X-X}(r,0)\setminus\{0\}.
\eea
It follows that $B$ is covered by $1+M^2(2r/\rho)^{2s}$ balls of
radius $\rho$.

Now suppose that $x\ne y$.  Let $a-b\in B\setminus\{0\}$. Again
$a\ne b$ and therefore $a$ is orthogonal to $b$. We have
$$
\|(a-b)-(x-y)\|^2=\{\begin{array}{l}\|a-x\|^2+\|b-y\|^2\cr
\|a+y\|^2+\|2x\|^2\cr \|2y\|^2+\|b+x\|^2\cr
\end{array}\quad\mbox{if}\quad\begin{array}{l}a\neq
y,\ b\neq x\cr a\neq y,\ b=x\cr a=y,\ b\neq x,\end{array}\right.
$$
and so
$$
\left.\begin{array}{ll}a\in B_X(r,x)&b\in B_X(r,y)\cr a\in
B_X(r,-y)&b\in B_X(r,x)\cr a\in B_X(r,y)&b\in B_X(r,-x)\cr a\in
B_X(r,y)&b\in
B_X(r,x)\end{array}\}\quad\mbox{if}\quad\{\begin{array}{l}a\neq y,\
b\neq x\cr a\neq y,\ b=x\cr a=y,\ b\neq x,\cr a=y,\
b=x.\end{array}\right.
$$
Therefore
\bea
    B\setminus\{0\}&\subseteq&\big(B_X(r,x)-B_X(r,y)\big)\cup
        \big(B_X(r,-y)-B_X(r,x)\big)\\
        &&\cup\big(B_X(r,y)-B_X(r,-x)\big)\cup
        \big(B_X(r,y)-B_X(r,x)\big).
\eea
Cover each of $B_X(r,x)$, $B_X(r,-x)$, $B_X(r,y)$ and $B_X(r,-y)$
by $M(2r/\rho)^s$ balls of radius $\rho/2$. An argument similar to
before yields a cover of $B$ by $1+4M^2(2r/\rho)^{2s}$ balls of
radius $r/2$.

Since we have $N_{X-X}(r,\rho)\le 1+4M^2(2r/\rho)^{2s}$ it follows
that $d_{\rm A}(X-X)\le 2s$.\end{proof}

\section{Non-existence of bi-Lipschitz linear embeddings}

In this section we give a simple example showing that if we
require a linear embedding (as in Theorem \ref{main}) then we can
do no better than almost bi-Lipschitz. First we prove the
following simple decomposition lemma for linear maps from $H$ onto
$\Re^k$ (cf.~comments in Hunt \& Kaloshin, 1997).

\begin{lemma}\label{decomp}
Suppose $L:H\rightarrow\Re^k$ is a linear map with
$L(H)=\Re^k$. Then $U=(\ker L)^\perp$ has dimension $k$, and $L$
can be decomposed uniquely as $MP$, where $P$ is the orthogonal
projection onto $U$ and $M:U\rightarrow\Re^k$ is an invertible
linear map.
\end{lemma}

Note that the result of this lemma shows Theorem
\ref{main} remains true with linear maps replaced by orthogonal
projections.  This gives a much more concise proof of the result in
Friz \& Robinson (1999).

\begin{proof}
Let $U=(\ker L)^\perp$ and suppose that there exist $m>k$ linearly
independent elements $\{x_j\}_{j=1}^m$ of $U$ for which
$Lx_j\neq0$. Then $\{Lx_j\}$ are elements of $\Re^k$; since $m>k$
at least one of the $\{Lx_j\}$ can be written as a linear
combination of the others:
$$
Lx_i=\sum_{j\neq i} c_j(Lx_j).
$$
It follows that
$$
\Big(x_i-\sum_{j\neq i}c_jx_j\Big)=0,
$$
which contradicts the definition of $U$.

Let $P$ denote the orthogonal projection onto $U$, and $M$ the
restriction of $L$ to $U$. Let $x\in H$, and decompose $x=u+v$,
where $u\in U$ and $v\in\ker L$.
Note that this decomposition is unique.
Clearly $Lx=Lu=Mu=M(Px)$. It remains to show that $M$ is invertible.
This is clear since $\dim U=\dim\Re^k=k$ and $M$
is linear.\end{proof}

Following Ben-Artzi et al.~(1993) we now prove

\begin{lemma}\label{noP}
Suppose that $X-X$ contains a set of the form
$\{\alpha_ne_n\}_{n=1}^\infty$ with $\{e_n\}_{n=1}^\infty$ an
orthonormal set. Then no linear map into any $\Re^k$ can be
bi-Lipschitz between $X$ and its image.
\end{lemma}

\begin{proof}
We assume that $L(H)=\Re^k$, otherwise it is possible to prune
some redundant dimensions from $\Re^k$. Suppose that $L$ is
bi-Lipschitz from $X$ into $\Re^k$. Write $L=MP$ as in Lemma
\ref{decomp}. Since $L$ is bi-Lipschitz on $X$ then for all $y\in
X-X$ we have
$$
\|y\|\le c|Ly|\quad=\quad c|MPy|\le C\|Py\|,
$$
where $C=c\|M\|$. In particular we have
$$
\|\alpha_n e_n\|\le c\|P(\alpha_ne_n)\|\qquad\Rightarrow\qquad
c\|Pe_n\|\ge 1.
$$
But
$$
k=\mbox{rank}\,P=\mbox{Trace}\,P\ge\sum_{n=1}^\infty(Pe_n,e_n)
    =\sum_{n=1}^\infty\|Pe_n\|^2=+\infty
$$
a contradiction.
\end{proof}

We note that this result also follows from Lemma 2.4 in
Movahedi-Lankarani \& Wells (2005) which gives a characterisation
of sets $X$ that can be linearly bi-Lipschitz embedded into some
$\Re^k$: such an embedding is possible if and only if the weak
closure of
$$
\{\frac{x-y}{\|x-y\|}:\ x,y\in X,\ x\neq y\}
$$
does not contain zero (``weak spherical compactness of $X$").

Now consider the homogeneous set $X=\{2^{-n}e_n\}\cup\{0\}$, which
has $d_{\rm A}(X)=0$. Since $X$ is an orthogonal sequence, it
follows that $X-X$ (which in particular contains $X$) is also
homogeneous; but Lemma \ref{noP} shows that no linear map into any
finite-dimensional Euclidean space can be bi-Lipschitz on $X$.
This shows that, with the requirement of linearly, our Theorem
\ref{main} cannot be improved.

However, note that there is a simple {\it nonlinear} bi-Lipschitz
map $\phi$ from $X$ into $[0,1]$, given by
$$
\phi(2^{-n}e_n)=2^{-n}:
$$
for $n<m$ we have
$$
\underbrace{\frac{1}{4}(2^{-n}+2^{-m})}_{\frac{1}{4}|2^{-n}e_n-2^{-m}e_m|}\le
2^{-(n+1)}\le\underbrace{|2^{-n}-2^{-m}|}_{|\phi(2^{-n}e_n)-\phi(2^{-m}e_m)|}\le
2^{-n}\le\underbrace{(2^{-n}+2^{-m})}_{|2^{-n}e_n-2^{-m}e_m|}.
$$
The relationship between linear embeddings and general
bi-Lipschitz embeddings is delicate. Suppose that $X$ is a
connected set containing more than one point. The result of
Proposition \ref{nasty} shows that even if $X$ can be linearly
bi-Lipschitz embedded into some $\Re^n$ it is nevertheless
bi-Lipschitz equivalent to a space $\phi(X)$ that cannot be
bi-Lipschitz embedded into any $\Re^n$ using a linear map.

\begin{blank}{

\section{Finite-dimensional results}

\begin{theorem}
Assume that $X$ is a compact subset of $\Re^M$ for some $M$ such
that $d_{\rm A}^\alpha(X-X)\le d$ for some $\alpha\ge0$. Then
provided that $n>(1+\alpha)d+1$, almost every linear map
$L:\Re^M\rightarrow\Re^n$ is injective on $X$ and there exists a
$\delta>0$ such that
$$
|x-y|\le c|Lx-Ly|(-\log|Lx-Ly|)\quad\mbox{for all}\quad x,y\in
X\mbox{ with }|x-y|\le\delta,
$$
where $c=2/\log 2$.
\end{theorem}

Without the assumption on the dimension of $A-A$ we can obtain the
following, reducing the standard embedding dimension by a factor
of $1/2$.

\begin{theorem}
Assume that $X$ is a compact subset of $\Re^M$ for some $M$ such
that $d_{\rm A}^\alpha(X)\le d$ for some $\alpha>0$. Then given
$n>d$ and $\theta<1-\frac{d}{n}$, almost every linear map
$L:\Re^M\rightarrow\Re^n$ is injective on $X$ and there exists a
$\delta>0$ such that
$$
|x-y|\le c|Lx-Ly|^\theta\quad\mbox{for all}\quad x,y\in X\mbox{
with }|x-y|\le\delta.
$$
\end{theorem}

}
\end{blank}

\section{Conclusion}

We have identified a new class of {\it almost homogeneous} metric
spaces, and shown that such spaces enjoy almost bi-Lipschitz
embeddings into Hilbert space. Furthermore we have shown that any
compact subset $X$ of a Banach space with $X-X$ almost homogeneous
can embedded into a finite-dimensional Euclidean space is an almost
bi-Lipschitz way, and used this to deduce the same for any compact
metric space $(X,d)$ with ${\mathscr F}(X)-{\mathscr F}(X)$ almost
homogeneous, where ${\mathscr F}:X\rightarrow L^\infty(X)$ is the
isometric Kuratowski embedding of $(X,d)$ into $L^\infty(X)$.

Some outstanding problems remain:
\begin{enumerate}
\item Is there a homogeneous subset of a Hilbert space that
cannot be bi-Lipschitz embedded into any $\Re^k$?
\item Can any (almost) homogeneous subset of a Hilbert
space be (almost) bi-Lipschitz embedded into some $\Re^k$?
\item Can one construct an almost bi-Lipschitz
embedding $f$ of a compact almost homogeneous metric space $(X,d)$
into a Hilbert space in such a way that $X-X$ is almost homogeneous?
(This would answer (2) positively.)
\item Is the exponent $\gamma$ in Theorem \ref{main} (the power of the $\slog$ term) in
any way optimal?
\item Can one bound the Assouad dimension of the
attractors of dissipative PDEs (or preferably the set of
differences of solutions lying on such attractors)?
\end{enumerate}

\begin{blank}{We remark here that a reduction of exponent $\gamma$ in
Theorem \ref{main} to less than $1/2$ along with the solution of (5)
would enable the construction of finite-dimensional ordinary
differential equations with unique solutions that accurately
reproduce the dynamics on the attractors of several significant
equations in mathematical physics.}\end{blank}

\bibliographystyle{amsplain}

\end{document}